\documentclass[11pt]{article}
\usepackage{graphicx} 

\title{\textbf{On slice measures of Green currents on $\cmplex\mathbb{P}^2$}}
\author{Christophe Dupont and Virgile Tapiero}
\date{19th december 2023}

\usepackage[T1]{fontenc}
\usepackage[english, french]{babel}
\usepackage{hyperref}
\usepackage{amsmath, amssymb, amsthm}
\usepackage[all]{xy} 
\usepackage{esint}
\usepackage{color}

\newtheorem{thm}{Theorem}[section]

\newtheorem{lemme}[thm]{Lemma}
\newtheorem{prop}[thm]{Proposition}

\newtheorem{defn}[thm]{Definition}

\newtheorem{rmq}[thm]{Remark}

\newtheorem{question}[thm]{Question}

\newtheorem*{theorem*}{Theorem}
\newtheorem*{theoremA*}{Theorem A}
\newtheorem*{theoremA.I*}{Theorem A.I}
\newtheorem*{theoremfrA.I*}{Théorème A.I}
\newtheorem*{theoremfrA.II*}{Théorème A.II}
\newtheorem*{theoremB*}{Theorem B}
\newtheorem*{theoremB.I*}{Theorem B.I}
\newtheorem*{theoremfrB.I*}{Théorème B.I}
\newtheorem*{theoremfrB.II*}{Théorème B.II}
\newtheorem*{theoremfrB.I.bis*}{Théorème B.I.bis}
\newtheorem*{propositionC*}{Proposition C}
\newtheorem*{propositionC.I*}{Proposition C.I}
\newtheorem*{corollaryD*}{Corollary D}
\newtheorem*{corollaryfrD.I*}{Corollaire D.I}
\newtheorem*{corollaryE*}{Corollary E}
\newtheorem*{corollaryfrE.I*}{Corollaire E.I}
\newtheorem*{corollaryC*}{Corollary C}

\newcommand{\reels}{\mathbb{R}}
\newcommand{\cmplex}{\mathbb{C}}
\newcommand{\LLog}{\ \mathrm{Log}}

\providecommand{\keywords}[1]
{
  {\small \noindent\textbf{\textit{Keywords---}} #1}
}

\usepackage[left=3.25cm,right=3.25cm,top=2cm,bottom=2cm]{geometry}

\begin{document}
\hypersetup{pdfborder=0 0 0}

\maketitle

        \begin{otherlanguage}{english}
            \begin{abstract}
           Let $f$ be a holomorphic map of $\cmplex\mathbb{P}^2$ of degree $d\geq 2$, let $T$ be its Green current and $\mu=T\wedge T$ be its equilibrium measure. We give a new proof of a theorem due to Dujardin asserting that $\mu\ll T\wedge\omega_{\mathbb{P}^2}$ implies $\lambda_2=\frac{1}{2}\LLog\ d$, where $\lambda_1 \geq \lambda_2$ are the Lyapunov exponents of $\mu$. Then, assuming $\mu\ll T\wedge\omega_{\mathbb{P}^2}$, we study slice measures $\nu :=T\wedge dd^c|W|^2$, where $W$ is a holomorphic local submersion. We give sufficient conditions on the Radon-Nikodym derivative of $\mu$ with respect to the trace measure $T\wedge\omega_{\mathbb{P}^2}$ ensuring $\mu=\nu$. The involved submersion $W$ comes from normal coordinates for the inverse branches of the iterates of $f$.
            \end{abstract}
        \end{otherlanguage}



\section{Introduction}

 Let $f$ be a holomorphic map of $\mathbb{P}^2$ of degree $d\geq2$. The Green current $T$ and the equilibrium measure  $\mu$ are invariant objects encoding the dynamical properties of $f$, we refer to the books  {\cite{dinsib10}}, {\cite{sib99}} by Dinh and Sibony. We recall that $T:=\lim_n\frac{1}{d^n}{f^n}^*\omega_{\mathbb{P}^2}$, where $\omega_{\mathbb{P}^2}$ is the normalized Fubini-Study $(1,1)-$form of $\mathbb{P}^2$, this is a closed positive current, with local H\"older potentials. It satisfies $f^*T=d T$ and its auto-intersection $\mu :=T\wedge T$ defines an invariant probability measure on  $\mathbb{P}^2$, which is mixing and  satisfies $f^*\mu=d^2 \mu$.\\

Berteloot-Loeb \cite{BL01} proved that $T$ is a smooth and non degenerate positive $(1,1)$-form on a non empty open subset of $\mathbb{P}^2$ if and only if $f$ is a Latt\`es map. Berteloot-Dupont \cite{berdup05} established later  that Berteloot-Loeb's condition on $T$ characterizes the condition $\mu\ll \mathrm{Leb}_{\mathbb{P}^2}$. Another characterization of  Latt\`es maps involves the Lyapunov exponents $\lambda_1\geq\lambda_2$ of $\mu$ (see Theorem \ref{thm:oseledec} for their definitions). Briend-Duval {\cite{BriDuv99}} proved that $\lambda_2\geq\frac{1}{2}\LLog\ d$. The equality $\lambda_1=\lambda_2=\frac{1}{2}\LLog\ d$ holds if and only if $\mu\ll \mathrm{Leb}_{\mathbb{P}^2}$, see \cite{dup06}. \\

It is now natural to characterize the mappings satisfying  $\lambda_1>\lambda_2=\frac{1}{2}\LLog\ d$. Hopefully, the presence of a minimal Lyapunov exponent should be equivalent to regular properties for $T$ and $\mu$, and perhaps (in some sense) to the existence of a one dimensional Latt\`es-like factor. The following result, due to Dujardin, nicely fits into this program. 

\begin{thm}[Dujardin {\cite[Theorem 3.6]{Duj12}}]\label{thm:Dujardin} If $\mu\ll  \sigma_T := T\wedge \omega_{\mathbb{P}^2}$ then $\lambda_2=\frac{1}{2}\LLog\ d$.
\end{thm}

He also asked the question of the reverse implication, a partial answer is provided in \cite{tap22}. The measure $\sigma_T$ is called the trace of $T$ and carries its mass. The proof of Dujardin's theorem is based on a construction of a $df-$invariant sub-bundle $\mathcal{T}\subset T\mathbb{P}^2$ (called Fatou directions) of rank $\geq 1$ that satisfies for $\sigma_T-$almost every $x$ and for any $\vec v \in \mathcal{T}_x \backslash\{0\}$ :
\begin{equation}\label{eq:bundleTofDujardin}
   \limsup_{n \to + \infty} \frac{1}{n} \LLog \,  ||d_x f^n (\vec v) ||\leq \frac{1}{2}\LLog\ d.
\end{equation}
Indeed, using $\mu\ll \sigma_T$ and Oseledec's Theorem \ref{thm:oseledec}, for $\mu-$almost every $x$, the $\limsup$ in (\ref{eq:bundleTofDujardin}) tends to $\lambda_1$ or $\lambda_2$, which implies in both cases $\lambda_2\leq\frac{1}{2}\LLog\ d$. Briend-Duval inequality $\lambda_2\geq \frac{1}{2}\LLog\ d$ then implies the equality. We note that these arguments only require that $\mu\ll \sigma_T$ on an open subset $U$ charged by $\mu$. \\

On a chart $U\subset\mathbb{P}^2$ equipped with holomorphic coordinates $(Z,W)$ we have: 
$$\sigma_T\asymp(T\wedge dd^c|Z|^2 + T\wedge dd^c|W|^2).$$ 
By construction $\mathrm{Supp}(\mu)\subset\mathrm{Supp}(\sigma_T)$, hence if $\mu(U)>0$ then $(T\wedge dd^c|Z|^2)(U)>0$ or $(T\wedge dd^c|W|^2)(U)>0$. Actually in this case $T\wedge dd^c|Z|^2$ and $T\wedge dd^c|W|^2$ both charge $U$. The idea is that if $T\wedge dd^c|Z|^2\equiv0$ on $U$, then the potentials of $T$ would be harmonic on almost every vertical disc contained in $U$, which implies that $\mu$ is null on $U$, see {\cite[\textsection 3.3]{dupont_rogue_2020}}.\\

The measures $T\wedge dd^c|Z|^2$ and $T\wedge dd^c|W|^2$ are called \textit{slices} of $T$. Dujardin's theorem gives an information on $\lambda_2$ when $\mu\ll\sigma_T$, the proof is obtained by applying general results concerning Fatou directions \cite{Duj12}. Our purpose in this article is to analyse more deeply the relations between $\mu$ and slices of $T$ when $\mu\ll\sigma_T$. \\

{\bf We first provide another proof of Dujardin's Theorem by using normal forms for the generic inverse branches of $f^n$ and forward recurrent properties.} We then provide another proof (of a weaker version of) Theorem \ref{thm:Dujardin} using backward recurrent properties (we indeed assume that the density of $\mu$ with respect to $T$ is bounded). The interest of this second proof is to introduce a decomposition of $\mu$ using normal coordinates. This decomposition is the cornerstone of the proof of Theorem \ref{thm:SlicePropertyINTRO} that we now introduce. \\

Relations between $\mu$ and slices of $T$ can be obtained when $f$ preserves a pencil of lines, given for instance by the meromorphic function $\pi[z:w:t]=[z:w]$. Dupont-Taflin {\cite{DT}} proved that in this situation $\mu$ and $T$ are related by the formula:
\begin{equation}\label{eq:formuledeDupont-Taflin}
    \mu = T\wedge \pi^{*}\mu_{\theta} \ \ \ \textrm{ (which implies } \mu \ll \sigma_T),
\end{equation}
where $\theta$ is the rational map satisfying $\pi\circ f= \theta\circ\pi$. We note that Jonsson  \cite{jon99} previously established an analogous formula for polynomial skew products on $\mathbb{C}^2$. On another hand, Berteloot-Loeb {\cite{BL}} proved that if $\theta$ is a Lattès map, then for every $a \in \mathbb{P}^1$ outside a finite subset, there exists a holomorphic coordinate $\zeta_a$ such that $\mu_{\theta}=dd^c|\zeta_a|^2$ on a neighborhood $V_a$ of $a$ ($\mu_{\theta}$ is the equilibrium measure of $\theta$). Combining this result with (\ref{eq:formuledeDupont-Taflin}), we obtain on $\pi^{-1}(V_a)$:
\begin{equation}\label{eq:eq:formuledeDupont-TaflinII}
    \mu = T\wedge dd^c|W|^2,\ \mathrm{with}\ W:=\zeta_a\circ\pi.
\end{equation}
This gives an example where a slice of $T$ is equal to $\mu$. Theorem \ref{thm:SlicePropertyINTRO} shows a similar formula. We assume $\mu\ll\sigma_T$ and denote $\psi\in L^1({\sigma_T})$ the Radon-Nikodym derivative of $\mu$ with respect to $\sigma_T$. Let $J:=\mathrm{Supp}(\mu)$ and $\psi |_J$ be the restriction of $\psi$ on $J$. Note that the measurable function $\psi$ is defined on $\mathbb{P}^2$, it is not unique since one can modify it outside $J$ by still preserving $\mu = \psi \sigma_T$. However, one has $\mu=\left(\psi\displaystyle{1\!\!1}_J\right)\sigma_T$ on $\mathbb{P}^2$ for every version of $\psi$. We denote 
$\psi_\Omega:=\psi \displaystyle{1\!\!1}_\Omega$ for every borel subset $\Omega$ of $\mathbb{P}^2$. We also denote by $\widehat x = (x_n)_{n \in \mathbb Z}$ the full orbits of the mapping $f$, by $\pi_0 : \widehat x \mapsto x_0$ the projection map and by $\widehat{\mu}$ the $\pi_0$-pullback of $\mu$ on the set of full orbits, see Section \ref{sec:normalforms}. In the next statement, the local submersion $W_{\widehat{x}}:B(x_0, \eta_{\varepsilon}(\widehat{x}))\to\mathbb C$ is provided by the normal form Theorem \ref{thm:normalforms}.

\begin{thm}\label{thm:SlicePropertyINTRO} Assume that $\mu=\psi_J  \sigma_T$  for some measurable function $\psi$ on $\mathbb{P}^2$ (in particular $\lambda_2 = \frac{1}{2}\LLog\ d$) and that $\lambda_1$ is strictly larger than $\lambda_2$. 

\begin{enumerate}
    \item If $\psi |_J$ is continuous, then for $\widehat{\mu}-$almost every $\widehat{x}$ there exist $r_1(\widehat{x})\leq\eta_{\varepsilon}(\widehat{x})$ and $C_{\widehat{x}} > 0$ such that:    
     $$\mu = C_{\widehat{x}} \left(T\wedge dd^c|W_{\widehat{x}}|^2\right)|_J\ \mathrm{on}\ B(x_0,r_1(\widehat{x})).$$
    \item If $\psi$ is continuous on an open neighborhood $V$ of $J$ and if $\mu=\psi_V\sigma_T$,
        then there exist $r_2(\widehat{x})\leq \eta_{\varepsilon}(\widehat{x})$ and $C_{\widehat{x}} > 0$ such that:
        $$\mu =  C_{\widehat{x}} \,  T\wedge dd^c|W_{\widehat{x}}|^2\ \mathrm{on}\ B(x_0,r_2(\widehat{x})).$$
\end{enumerate}
\end{thm}

The first Item implies that for $\widehat{\mu}-$almost every $\widehat{x}$, $\widehat{y}$ satisfying $\pi_0(\widehat{x}) = \pi_0(\widehat{y})$, the measures $\left(T\wedge dd^c|W_{\widehat{x}}|^2\right)|_J$ and $\left(T\wedge dd^c|W_{\widehat{y}}|^2\right)|_J$ coincide on a small ball centered at $\pi_0(\widehat{x})$. In particular, the germ of measure $\left(T\wedge dd^c|W_{\widehat{x}}|^2\right)|_J$ does not depend on the full orbit $\widehat{x}$ but only on the $\pi_0$-projection of $\widehat{x}$. A similar remark holds for the second Item.

Let us also observe that the hypothesis $\lambda_1 > \lambda_2$ in Theorem \ref{thm:SlicePropertyINTRO} is not restrictive. 
Assume indeed that $\mu=\psi_J  \sigma_T$ and that $\lambda_1= \lambda_2$. Then the Lyapunov exponents are both minimal (equal to $\frac{1}{2}\LLog\ d$) and $f$ is a Latt\`es mapping on $\mathbb{P}^2$. In this case there exists a finite ramified covering $\sigma : \mathbb{C}^2 / \Lambda \to \mathbb{P}^2$ such that $\sigma^*T$ is equal to a positive definite hermitian form $H$ with constant coefficients ($\Lambda$ being a cocompact real lattice in $\mathbb{C}^2$), see \cite{dup03}. This implies that, around every point $p$ outside an algebraic subset of $\mathbb{P}^2$ (the critical values of $\sigma$), the current $T$ coincides with $H$ in local coordinates given by $\sigma$. It follows that, for every local submersion $W$ defined near $p$, the measure $\mu = T \wedge T$ is equal to $T \wedge dd^c |W|^2$ multiplied by a smooth positive function. 

Let us finally note that the assumption of Item 2 morally replaces the function $\psi_J = \psi\displaystyle{1\!\!1}_J$ of the first Item (which is singular due to the multiplication by $\displaystyle{1\!\!1}_J$) by a continuous function on a neighborhood $V$ of $J$. More precisely, it requires the existence of an open neighborhood $V$ of $J$ such that $\psi|_V$ is continuous and such that $\mathrm{Supp}(\psi\sigma_T)|_V = J$. 
This assumption is fulfilled if one assumes that the restriction $\psi |_J$ is continuous and that $\mathrm{Supp} \,  T =J$, by applying Tietze-Urysohn theorem to $\psi|_J$. 

\begin{question}
    Is it sufficient to assume $\psi|_J$ continuous to get $\mu=T\wedge dd^c|W_{\widehat{x}}|^2$ locally ? 
\end{question}

We remark that the continuity of $\psi |_J$ is satisfied by suspensions of Latt\`es maps on $\mathbb{P}^1$, we provide explicit computations of $T$, $\mu$ and $\psi$ in Section \ref{sec:relevedeLattes}. \\

{\bf Acknowledgements :} The first author thanks the Simons Foundation, Laura DeMarco and Mattias Jonsson for the invitation at the Simons Symposium on Algebraic, Complex and Arithmetic Dynamics, in Schloss Elmau in August 2022.  His research is partially funded by the European Research Council (ERC GOAT 101053021). The two authors benefited from the support of the French government "Investissements d'Avenir" program integrated to France 2030 (ANR-11-LABX-0020-01). We thank the referee, whose advices and careful reading enabled us to improve the exposition of the article. 

\section{Classical results}

\subsection{Lyapunov exponents}

Let $f$ be a holomorphic mapping of $\mathbb{P}^2$ of degree $d\geq2$. Let $T$ be its Green current and $\mu =T\wedge T$ be its equilibrium measure, it is $f$-invariant and ergodic. Since $T$ has local continuous \emph{psh} potentials, $\mu$ integrates local \emph{psh} functions. In particular $\mu$ does not charge analytic subsets of $\mathbb{P}^2$ and $\mu(\mathcal C) =0$, where $\mathcal{C}:=\bigcup_{n\in\mathbb{Z}}f^n(\mathrm{Crit}\ f)$. It also allows to define the Lyapunov exponents $\lambda_1\geq\lambda_2$:
$$\lambda_1=\lim_{n\to+\infty}\frac{1}{n}\int_{\mathbb{P}^2}\LLog\ ||df^n||\ \mathrm{d}\mu \textrm{ and }  
\lambda_1+\lambda_2 = \int_{\mathbb{P}^2}\LLog\ \left|\mathrm{det}_{\cmplex}\ df\right|\ \mathrm{d}\mu.$$
Briend-Duval {\cite{BriDuv99}} proved that $\lambda_2\geq \frac{1}{2}\LLog\ d$. The Oseledec theorem states as follows. 

\begin{thm}[Oseledec]\label{thm:oseledec} \
\begin{enumerate}
    \item If $\lambda_1 = \lambda_2 = \lambda$, then there exists an invariant borel subset $A_{os}$ of full $\mu$-measure and disjoint from $\mathcal{C}$ such that for every $x \in A_{os}$:
    $$\forall \vec v \in T_x \mathbb{P}^2 \setminus \{ 0 \} \ , \ \lim_{n \to + \infty } \frac{1}{n} \LLog \, ||d_x f^n (\vec v)|| = \lambda.$$
    \item If $\lambda_1$ > $\lambda_2$, then there exists an invariant borel subset $A_{os}$ of full $\mu$-measure and disjoint from $\mathcal{C}$ such that for every $x \in A_{os}$, there exists $v_s(x) \in \mathbb{P}(T_x \mathbb{P}^2)$ satisfying:
$$\forall \vec v \in T_x \mathbb{P}^2 \setminus v_s(x) \ , \  \lim_{n \to + \infty} \frac{1}{n} \LLog \,  ||d_x f^n (\vec v)|| = \lambda_1.$$
$$\forall \vec v \in v_s(x) \setminus \{ 0 \}  \ , \ \lim_{n \to + \infty} \frac{1}{n} \LLog \, ||d_x f^n (\vec v)|| = \lambda_2.$$
Moreover $v_s$ is measurable and satisfies $[d_x f] (v_s(x)) = v_s(f(x))$ for every $x \in A_{os}$, where $[d_x f]$ is the projectivization of the tangent map $d_xf$. 
\end{enumerate}
\end{thm}

\subsection{Normal forms for inverse branches}\label{sec:normalforms}

We denote $X:=\mathbb{P}^2\backslash \mathcal{C}$ which is totally invariant by $f$, and we denote $\widehat{X}$ the set of all sequences $\widehat{x}=(x_n)_{n\in\mathbb{Z}}$ of elements of $X$ that are orbits under the action of $f$ $i.e.$ $x_{n+1}=f(x_n)$. Equipped with the natural projection $\pi_0:\widehat{X}\longrightarrow X$ defined by $\pi_0(\widehat{x})=x_0$, the dynamical system $(X,f,\mu)$ admits a natural extension $(\widehat{X},\widehat{f},\widehat{\mu})$, where $\widehat{f}$ is the left shift on $\widehat{X}$ and $\widehat{\mu}$ is the unique probability measure on $\widehat{X}$ invariant by $\widehat{f}$ such that $(\pi_0)_*\widehat{\mu} = \mu$. The measure $\widehat{\mu}$ is ergodic as $\mu$, see  {\cite[Chapter 10]{CFS12}}.

Given an orbit $\widehat{x}\in\widehat{X}$, there exists a family of inverse branches $(f^{-n}_{\widehat{x}})_{n\in\mathbb{N}}$ such that $f^{-n}_{\widehat{x}}$ is defined on $B(x_0, \eta_{\varepsilon}(\widehat{x}))$ for every $n \geq 0$, and such that $\mathrm{Lip}(f^{-n}_{\widehat{x}})\leq \beta_{\varepsilon}(\widehat{x})e^{-n(\lambda_2-\varepsilon)}$. The functions $\widehat{x}\mapsto \eta_{\varepsilon}(\widehat{x})$ and $\widehat{x}\mapsto \beta_{\varepsilon}(\widehat{x})$ are $\varepsilon-$tempered, \emph{i.e.} they satisfy :
$$ \forall n \in \mathbb Z \ , \ e^{-|n|\varepsilon}\varphi(\widehat{x})\leq \varphi (\widehat{f}^n(\widehat{x}) )\leq e^{+|n|\varepsilon}\varphi(\widehat{x}).$$
We refer to Briend-Duval {\cite{BriDuv99}} (see also {\cite{buzzi2003, dup06}}) for details about the construction of the inverse branches.  The following result provides normal forms for those mappings, see {\cite{bertelootdupont2019, bdm07, jonvar02}}. 

\begin{thm}\label{thm:normalforms}
For every $\varepsilon > 0$ small enough with respect to the Lyapunov exponents $\lambda_1\geq\lambda_2$, there exists an invariant borel subset $\Lambda_{\varepsilon}$ of full  $\widehat{\mu}$-measure, $\varepsilon$-tempered functions ${\rho}_\varepsilon : {\Lambda_{\varepsilon}} \to \left]0,1\right]$, $\beta_\varepsilon, L_\varepsilon, M_{\varepsilon} : {\Lambda_{\varepsilon}} \to [1,+\infty[$ and a function $N : {\Lambda_{\varepsilon}} \to \mathbb{N}$ satisfying the following properties. For every $\widehat{x} \in {\Lambda_{\varepsilon}}$, there exist injective holomorphic mappings 
$$\xi_{\widehat{x}} : B(x_0, \eta_{\varepsilon}(\widehat{x})) \longrightarrow \mathbb{D}^2(\rho_{\varepsilon}(\hat{x}))$$ satisfying :
\begin{enumerate}
    \item $\xi_{\widehat{x}}(x_0)=0$ and $d_{x_0} {\xi_{\widehat{x}}}(v_s(x_0))$ is the vertical axis in $\mathbb{C}^2$,
    \item {$\forall p,q \in B({x_0}, \eta_\varepsilon (\widehat{x}))$ , $\frac{1}{2}\mathrm{dist}(p,q) \leq ||\xi_{\widehat{x}}(p)-\xi_{\widehat{x}}(q)|| \leq \beta_{\varepsilon}(\widehat{x}) \,  \mathrm{dist}(p,q)$ ,}

\item the following diagram commutes for every  $n \geq N(\widehat{x})$: 
\begin{equation}\label{eq:diagramofnormalforms}
    \xymatrix{
      B({x_{-n}}, \eta_{\varepsilon}(\widehat{x}_{-n})) \ar[d]_{ \xi_{\widehat{x}_{-n}} } & &  B({x_0}, \eta_\varepsilon (\widehat{x})) \ar[ll]_{f^{-n}_{\widehat{x}}} \ar[d]^{\xi_{\widehat{x}}}  \\
      \mathbb{D}^2(\rho_{\varepsilon} (\widehat{x}_{-n}))  & & \mathbb{D}^2(\rho_{\varepsilon}({\widehat{x}}))  \ar[ll]^{R_{n,\widehat{x}}}
    }
\end{equation}
\end{enumerate}
The mappings $R_{n,\widehat{x}}$ have the following form depending on $(\lambda_1,\lambda_2)$: 
\begin{enumerate}
    \item[(i)] If $\lambda_1 = \lambda_2 = \lambda$, then $R_{n,\widehat{x}}$ is a linear and satisfies
    $$ e^{-n(\lambda + \varepsilon) }||(z,w)||  \leq ||R_{n,\widehat{x}}(z,w)|| \leq e^{-n(\lambda - \varepsilon)}||(z,w)||.$$
    \item[(ii)] If $\lambda_1 = k\lambda_2$ for some $k \geq 2$ (we say that $\lambda_1$ and $\lambda_2$ are resonant), then $R_{n,\widehat{x}} (z,w) = (\alpha_{n,\widehat{x}} z, \beta_{n,\widehat{x}} w) + ( \gamma_{n,\widehat{x}} w^k , 0)$.
    \item[(iii)]  If $\lambda_1 \not\in \{k \lambda_2, k \geq 1\}$, then $R_{n,\widehat{x}} (z,w)  = (\alpha_{n,\widehat{x}} z, \beta_{n,\widehat{x}} w)$.
\end{enumerate}
Moreover, in the cases $(ii)$ and $(iii)$, we have 
 $$ e^{- n (\lambda_1 + \varepsilon)} \leq |\alpha_{n, \widehat{x}}| \leq  e^{- n (\lambda_1 - \varepsilon) } \ \ , \ \  |\gamma_{n, \hat{x}}| \leq M_{\varepsilon}(\widehat{x})  e^{- n (\lambda_1 - \varepsilon)} $$
and 
$$ e^{-n(\lambda_2 + \varepsilon)} \leq |\beta_{n, \hat{x}}| \leq  e^{-n(\lambda_2 - \varepsilon)}. $$ \end{thm}

\begin{rmq}\label{rmq:commutativediagramforgerms} The diagram \eqref{eq:diagramofnormalforms} commutes for every fixed $n\geq0$ by reducing the radius of the balls. The property $n \geq N(\widehat{x})$ actually  ensures that $f^{-n}_{\widehat{x}} (B({x_0}, \eta_\varepsilon (\widehat{x}))) \subset   B({x_{-n}}, \eta_{\varepsilon}(\widehat{x}_{-n}))$. 
\end{rmq}

We shall need the following control of $(f^n)^*\omega_{\mathbb{P}^2}$ in Section \ref{sec:minimalexponent}.

\begin{lemme}\label{lemma:control(f^n)^*omegacorpsdutexte} For every $\widehat{x}\in\Lambda_{\varepsilon}$
and  $n \geq 0$, there exists $0 < r_n \leq \eta_{\varepsilon}(\widehat{x})$ such that 
$$(f^n)^*\omega_{\mathbb{P}^2}\geq 4 \beta_{\varepsilon}(\widehat{x})^{-2}e^{2n(\lambda_2-2\varepsilon)}  \omega_{\mathbb{P}^2} \textrm{ on } B(x_0,r_n) .$$
\end{lemme}

\noindent\underline{\textbf{Proof:}} We assume to simplify that $\lambda_1>\lambda_2$ and that $\lambda_1, \lambda_2$ are not resonant, the case of equality and the resonant case can be treated similarly (up to modifying the function $\beta_{\varepsilon}$). Let $\varepsilon >0$ be small enough such that $\lambda_1 - \lambda_2 > 2 \varepsilon$, ensuring $|\alpha_{n,\widehat{x}_n}|^{-1}\geq |\beta_{n,\widehat{x}_n}|^{-1}$. According to Remark \ref{rmq:commutativediagramforgerms}, there exists $r_n > 0$ and $0< r_n'\leq \eta_{\varepsilon}(\widehat{x}_n)$ such that $f^{-n}_{\widehat{x}_n} = \xi_{\widehat{x}}^{-1}\circ R_{n,\widehat{x}_{n}}\circ \xi_{\widehat{x}_n}$ on $B(x_n,r_n')$ and $B(x_0,r_n)\subset f^{-n}_{\widehat{x}_n}(B(x_n,r_n'))$. Since  we have on $B(x_0,r_n)$:
$$(f^{n})^{*}\omega_{\mathbb{P}^2}\geq ||(df^n)^{-1} ||^{-2}\omega_{\mathbb{P}^2} $$
and
$$||(df^n)^{-1}||^{-2} = ||df^{-n}_{\widehat{x}_n}||^{-2} \geq ||d\xi_{\widehat{x}}^{-1}||^{-2}||dR_{n,\widehat{x}_n}||^{-2} ||d\xi_{\widehat{x}_n}||^{-2},$$
 Theorem \ref{thm:normalforms} implies 
$$(f^{n})^{*}\omega_{\mathbb{P}^2}\geq 4||dR_{n,\widehat{x}_n}||^{-2}\beta_{\varepsilon}(\widehat{x}_n)^{-2}\omega_{\mathbb{P}^2} \textrm{ on } B(x_0, r_n). $$
Since $|\alpha_{n,\widehat{x}_n}|^{-1}\geq |\beta_{n,\widehat{x}_n}|^{-1}$ and $\beta_\varepsilon$ is $\varepsilon$-tempered, we deduce  on $B(x_0,r_n)$:
$$(f^{n})^{*}\omega_{\mathbb{P}^2}\geq 4 |\beta_{n,\widehat{x}_n}|^{-2} \beta_{\varepsilon}(\widehat{x})^{-2}e^{-2n\varepsilon}\omega_{\mathbb{P}^2}.$$
The conclusion follows from $|\beta_{n,\widehat{x}_n}|^{-2}\geq e^{2n(\lambda_2-\varepsilon)}$.\qed \\

We shall also need the following lemma in Section \ref{sec:aroundsliceconjecture}.

\begin{lemme}\label{lemma:auxiliarylemma}
Assume that $\mu=\psi_J  \sigma_T$  for some measurable function $\psi$ on $\mathbb{P}^2$. Let $\widehat{x}\in\Lambda_{\varepsilon}$ and $B:=B(x_0, \eta_{\varepsilon}(\widehat{x}))$. Let $A_n:=\{p\in B \cap J \, , \,  \psi \circ f^{-n}_{\widehat{x}}(p)>0\}$. 
\begin{enumerate}
\item For every $n \geq 0$,  the measures $\sigma_T$ and $(f^{-n}_{\widehat{x}})^* \sigma_T$ are equivalent on $B$, and $A_n = A_0$ modulo a borel subset of zero $\sigma_T-$measure. 
\item If $A \subset B \cap J$ satisfies $\mu(A)=0$, then $\sigma_T (A \cap A_0) = 0$.
In particular, for $a,b \geq0$, 
\begin{enumerate}
    \item[(i)] If $\psi\circ f^{-n}_{\widehat{x}} \leq b$ $\mu$-a.e. on $B \cap J$, then $\psi\circ f^{-n}_{\widehat{x}} \leq b$ $\sigma_T$-a.e. on $B \cap J$.
    \item[(ii)] If $\psi\circ f^{-n}_{\widehat{x}} \geq a$ $\mu$-a.e. on $B \cap J$, then $\psi\circ f^{-n}_{\widehat{x}}\geq a$ $\sigma_T$-a.e. on $B \cap J \cap A_0$.
\end{enumerate}
\end{enumerate}
 \end{lemme}

\noindent\textbf{\underline{Proof :}} Since  $f^{-n}_{\widehat{x}} : B \to f^{-n}_{\widehat{x}}(B)$ is a biholomorphism (and up to shrinking $B$), there exists $0 < \alpha < \beta$ such that 
$$ \alpha \, \omega_{\mathbb{P}^2} \leq (f^{-n}_{\widehat{x}})^*\omega_{\mathbb{P}^2} \leq \beta \,  \omega_{\mathbb{P}^2} \textrm{ on } B .$$
Using that $T$ is $f-$invariant, we deduce that $\sigma_T$ and $(f^{-n}_{\widehat{x}})^* \sigma_T$ are equivalent on $B$.
Now for every $p \in B \cap J$ and every $r$ small enough such that $B(p,r) \subset B$, 
$$\mu (B(p,r)) =  \int_{B(p,r)} \psi_J \sigma_T . $$    
But $\mu (B(p,r)) = d^{2n} \mu(f^{-n}_{\widehat{x}}(B(p,r)))$ and
$$  \mu(f^{-n}_{\widehat{x}}(B(p,r))) = \int_{B(p,r)} \psi_J \circ f^{-n}_{\widehat{x}} \  (f^{-n}_{\widehat{x}})^*\sigma_T  \simeq_{\alpha, \beta} d^{-n}  \int_{B(p,r)}  \psi_J \circ f^{-n}_{\widehat{x}}  \,  \sigma_T , $$
where $\simeq_{\alpha , \beta}$ means that the equality holds up to a multiplicative constant between $\alpha$ and $\beta$. We deduce by dividing by $\sigma_T(B(p,r))$ (which is positive since $p \in  J \subset \textrm{Supp} \, T$):
$$\fint_{B(p,r)} \psi_J \, \sigma_T  \simeq_{ \alpha, \beta} d^{-n} \fint_{B(p,r)}  \psi_J \circ f^{-n}_{\widehat{x}}  \,  \sigma_T . $$
By Lebesgue Theorem \cite[Theorem 2.12]{mattila} (taking limits when $r$ tends to zero), we get for $\sigma_T-$a.e. $p \in B \cap J$ that $\psi(p) > 0$ if and only if $ \psi \circ f^{-n}_{\widehat{x}}(p) > 0$. The second item comes from $0 = \mu(A)  = \int_{A \cap A_0}  \psi_J  \,  \sigma_T$, which gives $\sigma_T (A \cap A_0) = 0$ since $\psi_J$ is positive on $A \cap A_0$. To obtain $(i)$ and $(ii)$, we apply this fact to $\{p\in B \cap J \, , \,  \psi \circ f^{-n}_{\widehat{x}}(p) >b\}$ and to $\{p\in B \cap J \, , \,  \psi \circ f^{-n}_{\widehat{x}}(p) < a \}$ (which both have zero $\mu$-measure by assumption) and the fact that $A_n = A_0$ modulo a borel subset of zero $\sigma_T-$measure (provided by the first item).\qed

\section{A new proof of Dujardin's theorem}\label{sec:minimalexponent}

We provide another proof of Theorem \ref{thm:Dujardin}. We shall use Theorem \ref{thm:normalforms} and forward dynamics. We assume that $\mu \ll \sigma_T = T \wedge \omega_{\mathbb{P}^2}$.
This implies that there exists $\psi\in L^1({\sigma_T})$ such that $\mu=\psi_J\sigma_T$.  We recall that $\lambda_2\geq\frac{1}{2}\LLog\ d$ for every endomorphism of $\mathbb{P}^2$, it remains to prove the reverse inequality $\lambda_2\leq \frac{1}{2}\LLog\ d$ in our case of absolute continuity. For every $n \geq 0$, the set 
$$E_n:=\left\{x\in\mathbb{P}^2\ : \lim_{\rho\to0^+}\fint_{B(x,\rho)}\psi_J\circ f^{n}\ \sigma_T = \psi_J\circ f^n(x)\right\}$$
is a borel subset of full $\sigma_T-$measure (and so of full $\mu-$measure) by Lebesgue Theorem. We define the borel set 
$$E:=\left(\bigcap_{n\in\mathbb{N}}E_n\right)\backslash\left(\bigcup_{n\geq1}\mathrm{Crit}(f^n)\right),$$
which is also of full $\mu-$measure ($\mu$ does not charge analytic subsets).
Let us define 
$$F:=E\cap \left\{\frac{1}{\tau}\leq \psi_J\leq \tau \right\}\subset J,$$ where $\tau>0$ is chosen large enough to have $\mu(F)>0$. According to Poincar\'e recurrence theorem, for $\widehat{\mu}-$almost every $\widehat{x}\in \pi_0^{-1}(F)$, there exists a sub-sequence $\widehat{x}_{n_k}$ which satisfies $\widehat{x}_{n_k}\in \pi_0^{-1}(F)$. Let us fix $\varepsilon >0$ and apply the normal form Theorem \ref{thm:normalforms}, it yields a borel subset $\Lambda_\varepsilon$ of full $\widehat{\mu}$-measure of good orbits. Let us fix $\widehat{x}\in\pi_0^{-1}(F) \cap \Lambda_\varepsilon$ and a sub-sequence $(n_k)_k$ as before. Let us also fix an integer $k$. Since $x_0\in E$ avoids the critical set of $f^{n_k}$, there exists $\rho_k>0$ small enough such that $f^{n_k}$ is injective on the ball $B(x_0,\rho_k)$. Then we have for any $\rho\in]0,\rho_k]$:
\begin{equation}\label{eq:mu(A_rho)=d^-2n_kxmu(f^n_k(A_rho))}
    \mu(A_{\rho}) = d^{-2n_k}\mu(f^{n_k}(A_{\rho})),
\end{equation} 
where $A_{\rho}:=B(x_0,\rho)$. Up to taking a smaller $\rho_k$, we can assume that
$$f^{n_k}(A_{\rho})\subset B(x_{n_k},\eta_{\varepsilon}(\widehat{x}_{n_k}))\ \mathrm{and}\ A_{\rho}\subset  B(x_0,r_{n_k}),$$
where $r_{n_k}$ comes from Lemma \ref{lemma:control(f^n)^*omegacorpsdutexte}. 
 From $\mu=\psi_J\sigma_T$, Formula (\ref{eq:mu(A_rho)=d^-2n_kxmu(f^n_k(A_rho))}) yields for every $\rho\in]0,\rho_k]$:
\begin{align*}
    \int_{A_{\rho}}\psi_J\ \sigma_T & = d^{-2n_k}\int_{f^{n_k}(A_{\rho})}\psi_J\ \sigma_T\\
                                            & = d^{-2n_k}\int_{A_{\rho}}(\psi_J\circ f^{n_k}) \, d^{n_k} \, T\wedge (f^{n_k})^*\omega_{\mathbb{P}^2}
\end{align*}
where the second equality uses the invariance of the Green current $(f^{n_k})^{*}T=d^{n_k}T$ and the fact that $f^{n_k}$ is injective on $A_{\rho}$. Since $A_{\rho}\subset B(x_{0},r_{n_k})$, we can use the estimate given by Lemma \ref{lemma:control(f^n)^*omegacorpsdutexte} with ${n}={n_k}$ to deduce  
\begin{align*}
    \int_{A_{\rho}}\psi_J\ \sigma_T & \geq d^{-n_k}\int_{A_{\rho}}(\psi_J\circ f^{n_k}) 4 \beta_{\varepsilon}(\widehat{x})^{-2}e^{2n_k(\lambda_2-2\varepsilon)}\ T\wedge\omega_{\mathbb{P}^2}\\
                                            & \ \ \ \  = 4 d^{-n_k}\beta_{\varepsilon}(\widehat{x})^{-2}e^{2n_k(\lambda_2-2\varepsilon)}\int_{A_{\rho}}(\psi_J\circ f^{n_k})\ \sigma_T.
\end{align*}
We chose $x_0\in J$, thus $\sigma_T(A_{\rho})>0$ and we have :
$$\fint_{A_{\rho}}\psi_J\ \sigma_T \geq 4 d^{-n_k}\beta_{\varepsilon}(\widehat{x})^{-2}e^{2n_k(\lambda_2-2\varepsilon)}\fint_{A_{\rho}}(\psi_J\circ f^{n_k})\ \sigma_T.$$
But we also have $x_0\in E\subset E_0\cap E_{n_k}$, thus we can take limits when $\rho\to0^+$ to get 
$$\psi_J(x_0) \geq 4 d^{-n_k}\beta_{\varepsilon}(\widehat{x})^{-2}e^{2n_k(\lambda_2-2\varepsilon)}\times \psi_J(f^{n_k}(x_0)).$$ 
Since $x_0$ and $x_{n_k}=f^{n_k}(x_0)$ belongs to $F\subset \{\frac{1}{\tau}\leq \psi_J\leq \tau\}$, we deduce that:
$$\tau \geq \frac{4}{\tau}\beta_{\varepsilon}(\widehat{x})^{-2}e^{-2n_k(\frac{1}{2}\mathrm{Log}\ d - \lambda_2 + 2\varepsilon)}.$$
Since this is true for every integer $k$, we get $\frac{1}{2}\LLog\ d - \lambda_2 + 2\varepsilon \geq 0$. We obtain $\lambda_2\leq \frac{1}{2}\LLog\ d$ as desired, since $\varepsilon>0$ is arbitrary small.

\section{Decomposition of \texorpdfstring{$\mu$}{TEXT}}

\subsection{Normal coordinates and absolute continuity}

According to Theorem \ref{thm:normalforms}, there exists  an invariant subset of full $\widehat{\mu}-$measure $\Lambda_{\varepsilon}$ such that for any $\widehat{x}\in\Lambda_{\varepsilon}$, there exist inverse branches 
$$ f^{-n}_{\widehat{x}}:B(x_0, \eta_{\varepsilon}(\widehat{x}))\longrightarrow B(x_{-n}, \eta_{\varepsilon}(\widehat{x}_{-n})),\ n\geq N(\widehat{x}),$$
and holomorphic charts 
\begin{equation}\label{eq: IB}
\xi_{\widehat{x}} = (Z_{\widehat{x}} ,W_{\widehat{x}}) : B(x_0, \eta_{\varepsilon}(\widehat{x})) \longrightarrow \mathbb{D}^2(\rho_{\varepsilon}(\hat{x}))
\end{equation}
such that $\frac{1}{2}\leq||d\xi_{\widehat{x}}||\leq \beta_{\varepsilon}(\widehat{x})$ on $B(x_0, \eta_{\varepsilon}(\widehat{x}))$ and 
$$ f^{-n}_{\widehat{x}} = \xi_{\widehat{x}_{-n}}^{-1} \circ R_{n,\widehat{x}} \circ \xi_{\widehat{x}} $$
for any $n\geq N(\widehat{x})$, where $R_{n,\widehat{x}}$ is a polynomial map $(\mathbb C^2 , (Z_{\widehat{x}} ,W_{\widehat{x}})) \to  (\mathbb C^2 , (Z_{\widehat{x}_{-n}} ,W_{\widehat{x}_{-n}}))$. When $\lambda_1 > \lambda_2$ are not resonant,  $R_{n,\widehat{x}}$ is a linear diagonal map with eigenvalues  $|\alpha_{n,\widehat{x}}|\simeq e^{-n(\lambda_1\pm\varepsilon)}$ and $|\beta_{n,\widehat{x}}|\simeq e^{-n(\lambda_2\pm\varepsilon)}$.
 Let $(\omega_{\widehat{x},AB})_{A,B\in\{Z,W\}}$ be the smooth functions  defined by the formula
\begin{equation}\label{eq:expressionomegaintothenormalforms}
    \omega_{\mathbb{P}^2}\ = \sum_{A,B\in\{Z,W\}}\omega_{\widehat{x},AB} \, \frac{i}{2}\left(dA_{\widehat{x}}\wedge d\overline{B}_{\widehat{x}}\right)\ \mathrm{on}\ B(x_{0}, \eta_{\varepsilon}(\widehat{x})).
\end{equation}
For every $A, B$ in $\{Z,W\}$ and for every $n\geq N(\widehat{x})$ we define the following  functions:
$$ F_{\widehat{x},n}^{AB}  := \omega_{\widehat{x}_{-n},AB}\circ f^{-n}_{\widehat{x}}  : B(x_0,\eta_{\varepsilon}(\widehat{x}))\longrightarrow \cmplex . $$

\begin{prop}\label{lemma:lemmediese} Assume that $\mu \ll \sigma_T$ and let $\psi_J$ be the Radon-Nikodym derivative of $\mu$ with respect to $\sigma_T$. We assume  that the Lyapunov exponents $\lambda_1>\lambda_2$ are not resonant. For every $\widehat{x}\in\Lambda_{\varepsilon}$, for every borel set $U\subset B(x_0, \eta_{\varepsilon}(\widehat{x}))$ and for every $n\geq N(\widehat{x})$ we have:
\begin{align*}
    \mu(U) = d^{n}\int_{U}(\psi_J\circ f^{-n}_{\widehat{x}})\sum_{A,B\in\{Z,W\}} F_{\widehat{x},n}^{AB} \, \alpha^A_{n,\widehat{x}} \, \overline{\alpha^B_{n,\widehat{x}}}\ \ T\wedge \frac{i}{2}dA_{\widehat{x}}\wedge d\overline{B}_{\widehat{x}},
\end{align*}
where $\alpha_{n,\widehat{x}}^Z:=\alpha_{n,\widehat{x}}$ and $\alpha_{n,\widehat{x}}^W:=\beta_{n,\widehat{x}}$.
\end{prop}
\noindent\textbf{\underline{Proof:}}
From the equalities $\mu(U) = d^{2n}\mu(f^{-n}_{\widehat{x}}(U))$ and $\mu = \psi_J(T\wedge\omega_{\mathbb{P}^2})$, we obtain 
$$\mu(U) = d^{2n}\int_{U}(\psi_J\circ f^{-n}_{\widehat{x}})\  (f^{-n}_{\widehat{x}})^*T\wedge (f^{-n}_{\widehat{x}})^*\omega_{\mathbb{P}^2}.$$
By invariance of the Green current, we get $(f^{-n}_{\widehat{x}})^*T=d^{-n}T$ on $B(x_{0}, \eta_{\varepsilon}(\widehat{x}))$ and
$$\mu(U) = d^{n}\int_{U}(\psi_J\circ f^{-n}_{\widehat{x}})\  T\wedge (f^{-n}_{\widehat{x}})^*\omega_{\mathbb{P}^2}.$$
We conclude by pulling back Formula (\ref{eq:expressionomegaintothenormalforms}) ($\widehat{x}$ being replaced by $\widehat{x}_{-n}$) by $f^{-n}_{\widehat{x}}$ and using the fact that $R_{n,\widehat{x}}$ is linear and diagonal with eigenvalues $\alpha_{n,\widehat{x}}$ and $\beta_{n,\widehat{x}}$.\qed 

\begin{defn}\label{defn:decompositionofmu} Assume that $\mu\ll\sigma_T$ and let $\psi_J = \psi\displaystyle{1\!\!1}_J$ be the Radon-Nikodym derivative of $\mu$ with respect to $\sigma_T$. We assume  that the Lyapunov exponents  $\lambda_1>\lambda_2$ are not resonant. For every $U\subset B(x_0, \eta_{\varepsilon}(\widehat{x}))$ and  $n\geq0$, we have according to Proposition \ref{lemma:lemmediese}:
\begin{equation*}\label{eq:decompositionmu(U)=I_n(U)+J_n(U)}
    \mu(U) = I_n(U\cap J) + J_n(U\cap J),
\end{equation*}
where for every $D\subset B(x_0, \eta_{\varepsilon}(\widehat{x}))$:
\begin{equation*}\label{eq:expressiondeI_n(A)dansladecompositionmu=I_n+J_n}
    I_n(D) :=\ d^{n}|\beta_{n,\widehat{
        x}}|^2\int_{D}(\psi\circ f^{-n}_{\widehat{x}}) F_{\widehat{x},n}^{WW}  \  T \wedge \frac{i}{2}(dW_{\widehat{x}} \wedge d\overline{W}_{\widehat{x}})
\end{equation*}
and 
\begin{equation*}\label{eq:expressiondeJ_nquandlesexposantsneresonnentpas}
    J_n(D) :=\ d^{n} \sum_{(A,B)\ne(W,W)} \alpha^A_{n_k,\widehat{
        x}}\overline{\alpha^B_{n,\widehat{
        x}}}\int_{D}(\psi\circ f^{-n}_{\widehat{x}}) F_{\widehat{x},n}^{AB} \ T\wedge \frac{i}{2}(dA_{\widehat{x}}\wedge d\overline{B}_{\widehat{x}}) 
\end{equation*}
where $\alpha^Z_{n,\widehat{x}}:=\alpha_{n,\widehat{x}}$ and $\alpha^W_{n,\widehat{x}} := \beta_{n,\widehat{x}}$.
\end{defn}
\begin{rmq}\label{rk: def} We observe that if $\mu=\psi_V\sigma_T$ as in the second item of Theorem \ref{thm:SlicePropertyINTRO}, then
$$\mu(U) = I_n(U\cap V) + J_n(U\cap V)$$
for every $U\subset B(x_0,\eta_{\varepsilon}(\widehat{x}))$.
\end{rmq}

\subsection{Recurrence in \texorpdfstring{$\mathcal{A}_{\tau}$}{TEXT}}\label{sec:recurrenceinA_tau}

\begin{defn}\label{rouge} For $\tau>0$ let $\mathcal{A}_{\tau} \subset \Lambda_{\varepsilon}$ be the subset consisting of elements $\widehat y$ satisfying
\begin{enumerate}
\item $|\omega_{\widehat{y},AB} |\leq \tau$,
\item $\frac{1}{\tau}\leq\omega_{\widehat{y},AA} \leq \tau$, 
\item $S_{AB}(\widehat{y}) := \sup_{B(y_0,\eta_{\varepsilon}(\widehat{y}))}||d\omega_{\widehat{y},AB}||\leq\tau$,
\item $\psi_J \circ\pi_0 (\widehat y) \geq\frac{1}{\tau}$, 
\end{enumerate}
the two first items holding on $B(y_0, \eta_{\varepsilon}(\widehat{y}))$, the three first ones for every $A,B\in\{Z,W\}$. 
\end{defn}
Since  the distortion of the coordinates $\xi_{\widehat{x}} = (Z_{\widehat{x}} ,W_{\widehat{x}})$ is controlled by Theorem \ref{thm:normalforms}, there exists $\tau$ large enough such that $\mathcal{A}_{\tau}$ has positive $\widehat{\mu}$-measure. By Poincar\'e recurrence theorem, for $\widehat{\mu}-$almost every $\widehat{x}$, there exists $(n_k)_{k\in\mathbb{N}}$ such that $n_k\geq N(\widehat{x})$ and
$$\widehat{x}_{-n_k}\in \mathcal{A}_{\tau} .$$
We have specified $n_k\geq N(\widehat{x})$ to have $f^{-n_k}_{\widehat{x}}(B(x_0, \eta_{\varepsilon} (\widehat{x})) \subset B(x_{-n_k}, \eta_{\varepsilon}(\widehat{x}_{-n_k}))$. In this way the functions $F_{\widehat{x},n_k}^{AB}$ are well defined on $B(x_0, \eta_{\varepsilon}(\widehat{x}))$. By definition of $F_{\widehat{x},n}^{AB}$, we get from the first item of Definition \ref{rouge}:  
 \begin{equation}\label{eq:equationSIGMA}
            |  F_{\widehat{x},n_k}^{AB}  |\leq \tau \textrm{ on } B(x_0, \eta_{\varepsilon}(\widehat{x})).
 \end{equation}
Moreover, using the third item of Definition \ref{rouge} and $\mathrm{Lip}(f^{-n_k}_{\widehat{x}})\leq \beta_{\varepsilon}(\widehat{x})e^{-n_k(\lambda_2-\varepsilon)}$ (given by Theorem \ref{thm:normalforms}), we get
$$ ||d F_{\widehat{x},n_k}^{AB}||  \leq \tau \beta_{\varepsilon}(\widehat{x})e^{-n_k(\lambda_2-\varepsilon)} \textrm{ on } B(x_0, \eta_{\varepsilon}(\widehat{x})).$$
We deduce that 
\begin{equation}\label{lemma2:ascoliconvergence}
F_{\widehat{x},n_k}^{AB} \textrm{ uniformly converges on } B(x_0,\eta_{\varepsilon}(\widehat{x})) \textrm{ to a constant } F_{\widehat{x}}^{AB} \textrm{ up to a subsequence}
\end{equation}
and that (according to the second item of Definition \ref{rouge}):   
\begin{equation}\label{lemma:ascoliconvergence}
    F_{\widehat{x},n_k}^{AA} \geq \frac{1}{\tau} \textrm{ on } B(x_0,\eta_{\varepsilon}(\widehat{x})) . 
    \end{equation}
     
\subsection{A weaker version of Dujardin's theorem}
\label{sec:anotherproofofthefactthatabsolutecontinuityimpliesaminimalexponent}

We give another proof of Theorem \ref{thm:Dujardin} assuming that the Radon-Nikodym derivative $\psi_J\in L^{\infty}(\sigma_T)$. Our motivation is to use backward iterates, namely the inverse branches $f^{-n_k}_{\widehat{x}}$, in order to prepare the proof of Theorem \ref{thm:SlicePropertyINTRO} in Section $\ref{sec:proofundercontinuityassumption}$. As in Section \ref{sec:minimalexponent}, it suffices to prove that $\lambda_2\leq\frac{1}{2}\LLog\ d$.

For sake of simplicity, we assume that the Lyapunov exponents $\lambda_1 > \lambda_2$ are not resonant. This  implies that the polynomial map $R_{n,\widehat{x}}$ appearing in the diagram (\ref{eq:diagramofnormalforms}) is linear and diagonal. The proof can be adapted in the resonant case, see Section \ref{sec:abouttheresonantcase}. 

We apply the normal form Theorem \ref{thm:normalforms} with $\varepsilon >0$, it yields a borel subset $\Lambda_\epsilon$ of full $\hat \mu$-measure of good orbits.
We use the set $\mathcal{A}_{\tau}$  defined in Section \ref{sec:recurrenceinA_tau}. Let us fix a $\widehat{\mu}-$generic element $\widehat x \in \Lambda_\epsilon$, in particular $x_0 \in J$. Let $(n_k)_{k\in\mathbb{N}}$ be such that $\widehat{x}_{-n_k}\in\mathcal{A}_{\tau}$.\\

{\bf 1) Construction of a set of continuity $C$.} 

Let $\delta>0$ be such that $\mu(B(x_0,\eta_{\varepsilon}(\widehat{x}))) > 2\delta$. By Lusin's theorem there exists a compact set $C_0\subset B(x_0,\eta_{\varepsilon}(\widehat{x}))$ such that $\mu(B(x_0,\eta_{\varepsilon}(\widehat{x}))\backslash C_0) < \delta$ and $\psi_J$ is continuous on $C_0$. Using again Lusin's theorem we construct by induction a sequence $(C_n)_n$ of compact subsets of $B(x_0,\eta_{\varepsilon}(\widehat{x}))$ such that 
$$C_{n+1}\subset C_n\ \mathrm{and}\ \mu(C_n\backslash C_{n+1})<2^{-(n+1)}\delta$$
and
$$ \psi_J\circ f^{-(n+1)}_{\widehat{x}} \in C^0(C_{n+1}).$$
Then $\mu(C_n)>\mu(C_{n-1})-2^{-n}\delta>\cdots>\mu(C_0)-\delta>\mu(B(x_0,\eta_{\varepsilon}(\widehat{x})))-2\delta>0$. Thus $$C:= \bigcap_{n\in\mathbb{N}}C_n \subset B(x_0,\eta_{\varepsilon}(\widehat{x})) $$
satisfies $\mu(C)>0$. Let $\mu_{C}$ be the measure on $B(x_0,\eta_{\varepsilon}(\widehat{x}))$ defined by $$\mu_C (\cdot) := \mu(\cdot  \cap C) .$$
 Its support ${\rm supp} \,  \mu_C$ is included in $C \cap J$ because $C$ is a closed subset of $B(x_0,\eta_{\varepsilon}(\widehat{x}))$. \\

{\bf 2) Construction of $x_0'$ in ${\rm supp} \,  \mu_C$, with boundedness properties.} 

Let $G$ be a borel subset of full $\sigma_T-$measure (and so of full $\mu-$measure) satisfying
$$\forall p\in G,\ 0\leq\psi_J(p)\leq ||\psi_J||_{L^{\infty}(\sigma_T)} ,$$
and let $G_k$ be the borel set 
$$G_k:=f^{-n_k}_{\widehat{x}}(B(x_0,\eta_{\varepsilon}(\widehat{x})))\bigcap G.$$

Let us verify that $\mu(f^{n_k}(G_k)) = \mu(B(x_0,\eta_{\varepsilon}(\widehat{x})))$.
Since $\mu(G)=1$, we have $\mu(G_k) = \mu(f^{-n_k}_{\widehat{x}}(B(x_0,\eta_{\varepsilon}(\widehat{x})))) = d^{-2n_k} \mu(B(x_0,\eta_{\varepsilon}(\widehat{x})))$. But $f^{n_k}$ is injective on $f^{-n_k}_{\widehat{x}}(B(x_0,\eta_{\varepsilon}(\widehat{x})))$, which contains $G_k$, hence $\mu(f^{n_k}(G_k))=d^{2n_k}\mu(G_k)$, completing our verification.

 Therefore $$\mu(\bigcap_{k\in\mathbb{N}}f^{n_k}(G_k))=\mu(B(x_0,\eta_{\varepsilon}(\widehat{x}))) , $$
  which is positive since $x_0\in J$. Let 
\begin{equation}\label{eq:x'control}
     {x}'_0\in \bigcap_{k\in\mathbb{N}}f^{n_k}(G_k) \cap \left\{\frac{1}{\tau'}\leq \psi_J\right\}\cap {\rm supp} \,  \mu_C ,
\end{equation}
where $\tau'>1$ is large enough so that the set defined in (\ref{eq:x'control}) is not empty. Using the fact that 
$f^{-n_k}_{\widehat{x}} \circ (f^{n_k})_{\vert G_k} = {\rm Id}_{G_k} $, we observe that 
\begin{equation}\label{inG}
 \forall k\in\mathbb{N} \ , \ f^{-n_k}_{\widehat{x}}(x_0')\in G_k \subset G . 
 \end{equation}
 
 {\bf 3) Conclusion.} 

 For a given $k$, let $\rho_k>0$ be small enough so that for any $0 < \rho < \rho_k$:
$$A_{\rho}:=\left(B(x'_0,\rho) \cap C\right)\subset B(x_0,\eta_{\varepsilon}(\widehat{x})).$$
Since $||d\xi_{\widehat{x}}||\leq\beta_{\varepsilon}(\widehat{x})$ we can assume that $\sum_{A,B}|T\wedge \frac{i}{2}dA_{\widehat{x}}\wedge d\overline{B}_{\widehat{x}}|\leq \beta_{\varepsilon}(\widehat{x})^2\sigma_T$ on $B(x_0,\eta_{\varepsilon}(\widehat{x}))$. According to Proposition \ref{lemma:lemmediese}, $|F_{\widehat{x},n_k}^{AB}|\leq\tau$ (see Equation (\ref{eq:equationSIGMA})) and $|\alpha_{n,\widehat{x}}|\leq|\beta_{n,\widehat{x}}|$, we deduce
$$ \int_{A_{\rho}}\psi_J\ \mathrm{d}\sigma_T = \mu(A_{\rho}) \leq d^{n_k}|\beta_{n_k,\widehat{x}}|^2 \tau \beta_{\varepsilon}(\widehat{x}) ^2 \int_{A_{\rho}}(\psi_J\circ f^{-n_k}_{\widehat{x}})\ \mathrm{d}\sigma_T.$$
According to Equation (\ref{eq:x'control}), we have $x_0'\in {\rm supp} \,  \mu_C$, thus $\sigma_T(A_{\rho})>0$ and 
$$\fint_{A_{\rho}}\psi_J\ \mathrm{d}\sigma_T \leq d^{n_k}|\beta_{n_k,\widehat{x}}|^2  \tau \beta_{\varepsilon}(\widehat{x})^2\fint_{A_{\rho}}(\psi_J\circ f^{-n_k}_{\widehat{x}})\ \mathrm{d}\sigma_T.$$
Equation (\ref{eq:x'control}) also gives $\frac{1}{\tau'}\leq \psi_J(x'_0)$ and $x'_0 \in {\rm supp} \,  \mu_C \subset C\subset C_0\cap C_{n_k}$. We moreover have $A_{\rho} \subset C$ by definition. Hence, if $\rho$ tends to $0$, we get by continuity of $\psi_J$  (resp. $\psi_J\circ f^{-n_k}_{\widehat{x}}$) on $C_0$ (resp. on $C_{n_k}$):
$$\frac{1}{\tau'}\leq \psi_J(x'_0)\leq d^{n_k}|\beta_{n_k,\widehat{x}}|^2 \tau \beta_{\varepsilon}(\widehat{x})^2(\psi_J\circ f^{-n_k}_{\widehat{x}})(x'_0)\leq d^{n_k}|\beta_{n_k,\widehat{x}}|^2 \tau \beta_{\varepsilon}(\widehat{x})^2 ||\psi_J||_{L^{\infty}(\sigma_T)}, $$
the last inequality coming from the choice of $x_0'$, see Equation \eqref{inG}. We obtain for every $k$ :
$$\left( \tau' \tau \beta_{\varepsilon}(\widehat{x})^2 ||\psi_J||_{L^{\infty}(\sigma_T)}  \right)^{-1}\leq d^{n_k}|\beta_{n_k,\widehat{x}}|^2\leq e^{2n_k\left(\frac{1}{2}\mathrm{Log}\ d - \lambda_2 + \varepsilon\right)}.$$
Hence $\lambda_2\leq \frac{1}{2}\LLog\ d + \varepsilon$, which gives $\lambda_2\leq \frac{1}{2}\LLog\ d$ when $\varepsilon$ tends to zero, as desired.\qed 

\section{Estimates on \texorpdfstring{$I_n$}{TEXT} and \texorpdfstring{$J_n$}{TEXT} in the decomposition of \texorpdfstring{$\mu$}{TEXT}}\label{sec:aroundsliceconjecture}

\subsection{Study of \texorpdfstring{$I_n$}{TEXT} and \texorpdfstring{$J_n$}{TEXT}}\label{sec:decompositionofmuintothenormalforms}

We recall that \texorpdfstring{$I_n$}{TEXT} and \texorpdfstring{$J_n$}{TEXT} were introduced in Definition \ref{defn:decompositionofmu}.

\begin{prop}\label{lemma:decompositionofmuforthesliceconjecture} Assume that $\mu\ll \sigma_T$ and let $\psi\displaystyle{1\!\!1}_J$ be the Radon-Nikodym derivative of $\mu$ with respect to $\sigma_T$. Assume that $\lambda_1>\lambda_2$ are not resonant and let $\varepsilon > 0$ be small enough such that 
$$ \lambda_1 + \lambda_2 - 2 \varepsilon > \LLog \, d \ \ , \ \ \lambda_1 > \frac{1}{2}\ \mathrm{Log}\, d + \varepsilon . $$
 Let $\widehat{x}$ be a $\widehat \mu$-generic orbit and let $(n_k)_k$ be a sub-sequence such that $\widehat{x}_{-n_k}\in \mathcal{A}_{\tau}$. Let 
 $$\sigma_{\widehat{x}} := T\wedge dd^c|W_{\widehat{x}}|^2 . $$
  Let $B := B(x_0,\eta_\varepsilon(\widehat x))$, $A_0:=\{p\in B \cap J \, , \,  \psi (p)>0\}$ and $U \subset B$ be a borel set.
\begin{enumerate}
    
\item If  $\psi\circ f^{-n_k}_{\widehat{x}} \leq b \ \mu-a.e. \textrm{ on } B$ for every $k\geq0$, then there exists $C_{\widehat{x}}>0$ which does not depend on $U \subset B$ such that 
    $$ J_{n_k}(U \cap J ) \leq  C_{\widehat x} \, e^{-n_k (\lambda_1 + \lambda_2 - \mathrm{Log} \, d)} e^{2 n_k \varepsilon} \underset{k \to \infty}{\longrightarrow} 0.$$

\item If  $\psi\circ f^{-n_k}_{\widehat{x}} \leq b \ \mu-a.e. \textrm{ on } B$ for every $k\geq0$ and if $\mu(U)>0$, then 
$$\sigma_{\widehat{x}}(U\cap J \cap A_0) >0 .$$
    
\item If $0< a \leq \psi\circ f^{-n_k}_{\widehat{x}}\leq b \ {\mu}-a.e.\ \mathrm{on}\ B$ for every $k\geq0$, then
        $$d^{n_k}|\beta_{n_k,\widehat{x}}|^2 a \tau^{-1} \sigma_{\widehat{x}}(U\cap J \cap A_0) \leq I_{n_k}(U\cap J)  \leq d^{n_{k}}|\beta_{n_{k},\widehat{x}}|^{2}b  \tau \sigma_{\widehat{x}}(U\cap J).$$
        In particular $(d^{n_{k}}|\beta_{n_k,\widehat{x}}|^2)_k$ converges, up to a sub-sequence, to some  $u>0$. 
\end{enumerate}
\end{prop}

\begin{rmq} We assume in Proposition \ref{lemma:decompositionofmuforthesliceconjecture} that $\lambda_1 , \lambda_2$ are not resonant. The same statement actually holds in the resonant case, this is explained in Section \ref{sec:abouttheresonantcase}.
\end{rmq}

\noindent\textbf{\underline{Proof of Proposition \ref{lemma:decompositionofmuforthesliceconjecture} :}} \  Let us prove the first item. By Definition \ref{defn:decompositionofmu},
$$   J_{n_k}(U \cap J) =  \sum_{(A,B)\ne(W,W)} d^{n_k} \, \alpha_{n_k,\widehat{x}}^A \,  \overline{\alpha_{n_k,\widehat{x}}^B} \int_{U \cap J} (\psi\circ f_{\widehat{x}}^{-n_k})   F_{\widehat{x},n_k}^{AB} \ \ T\wedge dA_{\widehat{x}}\wedge d\overline{B}_{\widehat{x}} .$$
Let us introduce the measures (which are $\ll \sigma_T$): 
$$\sigma_{\widehat{x}} := T\wedge dd^c|W_{\widehat{x}}|^2 \ \ , \ \ \lambda_{\widehat{x}} := T\wedge dd^c|Z_{\widehat{x}}|^2 . $$ 
Using $|\alpha_{n,\widehat{x}}|\leq|\beta_{n,\widehat{x}}|$, the (assumed) upper bound on $\psi\circ f^{-n_k}_{\widehat{x}}$,  the upper bound on $F_{\widehat{x},n_k}^{AB}$ provided by (\ref{eq:equationSIGMA}), the item $2.(i)$ of Lemma \ref{lemma:auxiliarylemma} and Cauchy-Schwarz inequality, we get$$  |J_{n_k}(U \cap J)| \leq d^{n_k}|\alpha_{n_k,\widehat{x}}\beta_{n_k,\widehat{x}}| \,  b \, \tau \, (\sigma_{\widehat{x}}+\lambda_{\widehat{x}})(U  \cap J).$$
The first item then follows from $|\alpha_{n, \widehat{x}}| \leq  e^{- n (\lambda_1 - \varepsilon) }$ and $ |\beta_{n, \hat{x}}| \leq  e^{-n(\lambda_2 - \varepsilon)}$. \\

We show the second item. By using Item 1 of Lemma \ref{lemma:auxiliarylemma} (for the last equality), we get
$$\mu(U) = \mu(U\cap J) = d^{2n}\mu(f^{-n}_{\widehat{x}}(U\cap J)) =  d^{2n}\int_{U\cap J \cap A_0}(\psi\circ f^{-n}_{\widehat{x}}) \ (f^{-n}_{\widehat{x}})^*(T\wedge\omega_{\mathbb{P}^2}) . $$
Cauchy-Schwarz inequality implies (up to multiplication by a constant):
$$T\wedge\omega_{\mathbb{P}^2}\leq  \omega_{\widehat{x}_{-n},WW} \, \sigma_{\widehat{x}_{-n}} + \omega_{\widehat{x}_{-n},ZZ} \, \lambda_{\widehat{x}_{-n}} \ \mathrm{on}\ B(\widehat{x}_{-n},\eta_{\varepsilon}(\widehat{x}_{-n})) .$$
We now fix $n=n_k$. Using the (assumed) upper bound on $\psi\circ f^{-n_k}_{\widehat{x}}$ on $B$ combined with item $2.(i)$ of Lemma \ref{lemma:auxiliarylemma}, and the upper bounds on $F_{\widehat{x},n_k}^{WW} = \omega_{\widehat{x}_{-n},WW} \circ f^{-n_k}_{\widehat{x}}$, $F_{\widehat{x},n_k}^{ZZ} = \omega_{\widehat{x}_{-n},ZZ} \circ f^{-n_k}_{\widehat{x}}$ provided by Equation 
(\ref{eq:equationSIGMA}), we get
$$\mu(U)  \leq   b d^{2n_k} \tau \left( \int_{U\cap J \cap A_0}  (f^{-n_k}_{\widehat{x}})^*(\sigma_{\widehat{x}_{-n_k}}) + \int_{U\cap J \cap A_0} \ (f^{-n_k}_{\widehat{x}})^*(\lambda_{\widehat{x}_{-n_k}}) \right) . $$
We deduce $$\mu(U)  \leq  b d^{2n_k} \tau \left( \int_{U\cap J \cap A_0}d^{-n_k}|\beta_{n_k,\widehat{x}}|^{2}\ \sigma_{\widehat{x}}+    \int_{U\cap J \cap A_0}d^{-n_k}|\alpha_{n_k,\widehat{x}}|^{2}\ \lambda_{\widehat{x}} \right) , $$
hence
$$0<\mu(U) \leq  b d^{n_k} \tau \left(  |\beta_{n_k,\widehat{x}}|^2 \sigma_{\widehat{x}}(U\cap J \cap A_0)   +   |\alpha_{n_k,\widehat{x}}|^2 \lambda_{\widehat{x}}(U\cap J \cap A_0 ) \right). $$
To conclude we observe that 
$$d^{n_k} |\alpha_{n_k,\widehat{x}}|^2 \leq e^{n_k(\mathrm{Log}\ d\ -2\lambda_1\ +2\varepsilon)},$$
which tends to $0$ when $k$ tends to infinity. Hence $\sigma_{\widehat{x}}(U\cap J \cap A_0)>0$. \\

Let us prove the third item. According to item $2.(ii)$ of Lemma \ref{lemma:auxiliarylemma}, we have
$$\psi\circ f^{-n_k}_{\widehat{x}} \geq a,\ \sigma_{\widehat{x}}-a.e.\ \textrm{on}\ B \cap J \cap A_0 .$$ 
Hence, according to Definition \ref{defn:decompositionofmu}, we get
$$I_{n_k}(U\cap J \cap A_0)\geq d^{n_k}|\beta_{n_k,\widehat{x}}|^2 a \int_{U\cap J\cap A_0} F_{\widehat{x},n_k}^{WW} \ \sigma_{\widehat{x}} \geq d^{n_k}|\beta_{n_k,\widehat{x}}|^2 a \frac{1}{\tau} \sigma_{\widehat{x}}(U\cap J \cap A_0),$$
where the last inequality comes from Equation \eqref{lemma:ascoliconvergence}.  Similarly, by using item $2.(i)$ of Lemma \ref{lemma:auxiliarylemma} and Equation \eqref{eq:equationSIGMA}, we get 
$$  I_{n_k}(U\cap J) \leq d^{n_{k}}|\beta_{n_{k},\widehat{x}}|^{2}b  \tau \sigma_{\widehat{x}}(U\cap J) $$
This proves the two stated inequalities on $I_{n_k}(U\cap J)$. Let us verify the last fact. Since $0<\mu(B) = I_{n_k}(B\cap J) +  J_{n_k}(B\cap J)$ and since $J_{n_k}(B)$ tends to $0$ by the first item, we get $\frac{1}{2}\mu(B)\leq I_{n_k}(B\cap J)\leq 2\mu(B)$ for $k$ large enough, which implies 
$$0<\frac{1}{2b\tau} \frac{\mu(B)}{\sigma_{\widehat{x}}(B\cap J)}\leq d^{n_{k}}|\beta_{n_{k},\widehat{x}}|^2\leq \frac{2\tau}{a} \frac{\mu(B)}{\sigma_{\widehat{x}}(B\cap J \cap A_0)}<+\infty,$$
since $\sigma_{\widehat{x}}(B\cap J \cap A_0) > 0$ by the second item.\qed

\subsection{Application: an equivalence between \texorpdfstring{$\mu$}{TEXT} and a slice of \texorpdfstring{$T$}{TEXT} }\label{sec:proofunderboundnessassumption}

We prove that $\mu$ can be approximated by a slice of $T$ under absolute continuity assumptions. The proof is a direct application of the decomposition of $\mu$ fixed in Definition \ref{defn:decompositionofmu} and of Proposition \ref{lemma:decompositionofmuforthesliceconjecture}. 

\begin{prop}\label{thm:resultsofSliceProperty} Assume that $\mu\ll \sigma_T$ and  $\lambda_1>\lambda_2=\frac{1}{2}\LLog\ d$. Let $\psi\displaystyle{1\!\!1}_J$ be the Radon-Nikodym derivative of $\mu$ with respect to $\sigma_T$. Let $\widehat{x}$ be a $\widehat \mu$-generic element and let $(n_k)_k$ be a sub-sequence such that $\widehat{x}_{-n_k}\in \mathcal{A}_{\tau}$. Let $B := B(x_0,\eta_\varepsilon(\widehat x))$ and $A_0:=\{p\in B \cap J \, , \,  \psi (p)>0\}$.  If 
$$ 0< a \leq \psi\circ f^{-n_k}_{\widehat{x}}\leq b\ {\mu}-a.e.\ \mathrm{on}\ B \textrm{ for every } k \geq 0 , $$
then there exists $C_{\widehat{x}}>0$ such that on $B$:
$$C_{\widehat{x}}^{-1} (T\wedge dd^c|W_{\widehat{x}}|^2)|_{J \cap A_0}\leq  \mu \leq  C_{\widehat{x}} (T\wedge dd^c|W_{\widehat{x}}|^2)|_{J}.$$
\end{prop}

\noindent\textbf{\underline{Proof :}} Let $U\subset B$ be a borel set and let $I_k:=I_{n_k}(U\cap J)$, $J_k:=J_{n_k}(U\cap J)$. According to Definition \ref{defn:decompositionofmu}, we have $ \mu(U) = I_k + J_k$.
The first item of Proposition \ref{lemma:decompositionofmuforthesliceconjecture} yields $J_k\to 0$. Using the third item of Proposition  \ref{lemma:decompositionofmuforthesliceconjecture} (which in particular provides $\lim_k d^{n_k}|\beta_{n_k,\widehat{x}}|^2 = u>0$), Equation (\ref{eq:equationSIGMA}) and item $2.(i)$ of Lemma \ref{lemma:auxiliarylemma}, we get 
$$\mu(U)\leq u\, b \, \tau \, \sigma_{\widehat{x}}(U \cap J).$$
Similarly, by using Equation \eqref{lemma:ascoliconvergence} and item $2.(ii)$ of Lemma \ref{lemma:auxiliarylemma}, we get 
$$\mu(U)\geq\frac{u \, a}{\tau} \, \sigma_{\widehat{x}}(U\cap J \cap A_0).$$
We conclude by setting $C_{\widehat{x}}:=\max\{\tau/ua, ub\tau\}$. \qed 

\section{Proof of Theorem {\ref{thm:SlicePropertyINTRO}}}\label{sec:equalityoofmuwithslicesofT}

\subsection{Proof of the first item}\label{sec:proofundercontinuityassumption}

\begin{lemme}\label{prop:preparationdutheoremWeakSliceProperty} Assume that $\mu\ll \sigma_T$ and  $\lambda_1>\lambda_2=\frac{1}{2}\LLog\ d$.
Let $\psi\displaystyle{1\!\!1}_J$ be the Radon-Nikodym derivative of $\mu$ with respect to $\sigma_T$, and assume that $\psi|_J$ is continuous. 
 Let $\widehat{x}$ be a $\widehat \mu$-generic element and let $(n_k)_k$ be a sub-sequence such that $\widehat{x}_{-n_k}\in \mathcal{A}_{\tau}$ with $n_k\geq N(\widehat{x})$. There exists $r_1(\widehat{x}) \leq \eta_ {\varepsilon}(\widehat{x})$ and $x_{-\infty}\in J$ such that, up to a sub-sequence :
\begin{enumerate}
             \item $\lim_k d^{n_k}|\beta_{n_k,\widehat{x}}|^2 = u>0$ and $\lim_k J_{n_k}(U \cap J)= 0$ for any borel set $U \subset B(x_0,r_1(\widehat{x}))$,
             \item $\psi\circ f^{-n_{k}}_{\widehat{x}}$ uniformly converges to $\psi(x_{-\infty})$ on $B(x_0,r_1(\widehat{x}))\cap J$.
\end{enumerate}
\end{lemme}

\noindent\textbf{\underline{Proof:}} 
To prove the first item, in view of Proposition \ref{lemma:decompositionofmuforthesliceconjecture}, it suffices to verify 
$$\frac{1}{2\tau}\leq \psi\circ f^{-n_k}_{\widehat{x}} $$
on some ${B}(x_{0},r_1(\widehat{x})) \cap J$. By uniform continuity, there exists $\delta(\tau)>0$ such that 
$$     \forall a,b\in J,\ \mathrm{dist}(a,b)\leq\delta(\tau) \Longrightarrow |\psi(a)-\psi(b)| \leq \frac{1}{2\tau}. $$
Defining $r_1(\widehat{x}):=\min \left\{\eta_{\varepsilon}(\widehat{x}),\frac{\delta(\tau)}{\beta_{\varepsilon}(\widehat{x})} \right\}$, we have for every $p \in B(x_{0},r_1(\widehat{x}))$ :
$$\mathrm{dist}\left(f^{-n_k}_{\widehat{x}}(p) , f^{-n_k}_{\widehat{x}}(x_{0})\right)  \leq \mathrm{Lip}(f^{-n_k}_{\widehat{x}}) r_1(\widehat{x}) \leq \beta_{\varepsilon}(\widehat{x}) r_1(\widehat{x}) \leq \delta(\tau) . $$
Using $\psi(x_{-n_k})\geq \frac{1}{\tau}$ (coming from $\widehat{x}_{-n_k}\in\mathcal{A}_{\tau}$), we get $\psi\circ f^{-n_k}_{\widehat{x}} \geq \frac{1}{2\tau}$ on $B(x_0,r_1(\widehat{x})) \cap J$.

For the second item, the sequence $(x_{-n_k})_k \in J$ converges by compactness to some $x_{-\infty}\in J$. Since $\mathrm{Lip} \, f^{-n_{k}}_{\widehat{x}}$ tends to $0$, $\left(f^{-n_{k}}_{\widehat{x}}\right)_k$ converges uniformly on $B(x_0, r_1(\widehat{x}))$ to the constant mapping $x_{-\infty}$. The conclusion follows from the continuity of $\psi$ on $J$. 
\qed

$ $

We now prove the first item of Theorem \ref{thm:SlicePropertyINTRO}. Let $U\subset B(x_0,r_1(\widehat{x}))$ be a borel set. We recall that $\sigma_{\widehat{x}} =T\wedge dd^c|W_{\widehat{x}}|^2$ and that by Definition \ref{defn:decompositionofmu}:
$$\mu(U) = \mu(U\cap J) = d^{n_{k}}|\beta_{n_{k},\widehat{x}}|^2\int_{U\cap J}\left(\psi\circ f_{\widehat{x}}^{-n_{k}}\right)F_{\widehat{x},n_k}^{WW}\ \sigma_{\widehat{x}} + J_{n_k}(U\cap J) .$$
According to Equation \eqref{lemma2:ascoliconvergence}, $F_{\widehat{x},n_k}^{WW}$ uniformly converges on $B(x_0,r_1(\widehat{x}))$ to some constant $F_{\widehat{x}}$, up to a sub-sequence.
We conclude by using the two items of Lemma \ref{prop:preparationdutheoremWeakSliceProperty}, which yield $\mu(U)=C_{\widehat{x}} \, \sigma_{\widehat{x}}(U\cap J)$, where $C_{\widehat{x}}:=u  \, F_{\widehat{x}} \,  \psi({{x}_{-\infty}})$. Note that $C_{\widehat{x}}$ does not depend on $U$ and is positive by taking $U = B(x_0,r_1(\widehat{x}))$.

\subsection{Proof of the second item}\label{sec:ProofoftheSlicePropertywithcontinuityonaneighborhoodofSupp(mu)}

The following lemma is the counterpart of Lemma \ref{prop:preparationdutheoremWeakSliceProperty}, the function $\psi$ being now continuous on a neighborhood $V$ of $J$. The two items hold without any restriction to $J$. The proof is similar using Remark \ref{rk: def} (decomposition of $\mu$ on $V$), the uniform continuity of $\psi$ on $V$, and by introducing a radius $r_2(\widehat{x})$ satisfying $B(x_0,r_2(\widehat{x})) \subset V$. Note that $f^{-n_k}_{\widehat{x}}(B(x_0,r_2(\widehat{x})))\subset V$ for $k$ large enough.

\begin{lemme}\label{prop:preparationdutheoremWeakSlicePropertyII} Assume that $\mu\ll \sigma_T$ and $\lambda_1>\lambda_2=\frac{1}{2}\LLog\ d$. 
We assume that there exists a neighborhood $V$ of $J$ and a function $\psi\in L^1(\sigma_T)$ such that $\psi$ is continuous on $V$ and satifies $\mu=\psi_V \sigma_T$ (we recall that $\psi_V=\psi \displaystyle{1\!\!1}_V$). Let $\widehat{x}$ be a $\widehat \mu-$generic element and let $(n_k)_k$ be a sub-sequence such that $\widehat{x}_{-n_k}\in \mathcal{A}_{\tau}$ with $n\geq N(\widehat{x})$. There exists $r_2(\widehat{x}) \leq \eta_ {\varepsilon}(\widehat{x})$ (small enough to have $B(x_0,r_2(\widehat{x})) \subset V$) and $x_{-\infty}\in J$ such that, up to a sub-sequence :
\begin{enumerate}
             \item $\lim_k d^{n_k}|\beta_{n_k,\widehat{x}}|^2 = u>0$ and $\lim_k J_{n_k}(U)= 0$ for any borel set $U\subset B(x_0,r_2(\widehat{x}))$,
             \item $\psi\circ f^{-n_{k}}_{\widehat{x}}$ uniformly converges to $\psi(x_{-\infty})$ on $B(x_0,r_2(\widehat{x}))$.
\end{enumerate}
\end{lemme}

The second item of Theorem \ref{thm:SlicePropertyINTRO}  follows as in the end of Section \ref{sec:proofundercontinuityassumption}, by using
$$\mu(U) = d^{n_{k}}|\beta_{n_{k},\widehat{x}}|^2\int_{U}\left(\psi\circ f_{\widehat{x}}^{-n_{k}}\right)F_{\widehat{x},n_k}^{WW}\ \sigma_{\widehat{x}} + J_{n_k}(U) $$
for every $U \subset B(x_0,r_2(\widehat{x})) \subset V$, this formula is provided by Definition \ref{defn:decompositionofmu} and Remark \ref{rk: def}.

\section{About the resonant case}\label{sec:abouttheresonantcase}

We assume that $\lambda_1, \lambda_2$ are resonant: $\lambda_1=q\lambda_2$ for some $q \geq 2$. We assume that $\mu\ll \sigma_T$, let $\psi$ denote the Radon-Nikodym derivative of $\mu$ with respect to $\sigma_T$. Let us first explain how Proposition \ref{lemma:lemmediese} is affected. We shall compute the corresponding $J_n(D)$ appearing in Definition \ref{defn:decompositionofmu}.  According to Theorem \ref{thm:normalforms}, the map $R_{n,\widehat{x}}$ is equal to:
$$R_{n,\widehat{x}} (z,w) = (\alpha_{n,\widehat{x}} z, \beta_{n,\widehat{x}} w) + ( \gamma_{n,\widehat{x}} w^{q} , 0),$$
where $|\gamma_{n, \hat{x}}| \leq M_{\varepsilon}(\widehat{x})  e^{- n (\lambda_1 - \varepsilon)}$. 
Let $U\subset B(x_0,\eta_{\varepsilon}(\widehat{x}))$ be a borel set. Denoting $\psi_{n}:=\psi\circ f^{-n}_{\widehat{x}}$, we have $\mu(U) = I_n(U\cap J) + J_n(U\cap J)$, where $I_n(D)$ has the same expression as in the non resonant case :
$$I_n(D) =\ d^{n}|\beta_{n,\widehat{x}}|^2\int_{D} \psi_n F_{\widehat{x},n}^{WW}\  T\wedge \frac{i}{2}(dW_{\widehat{x}}\wedge d\overline{W}_{\widehat{x}}) . $$
In the resonant case, the formula giving $J_n(D)$ is modified as follows: 
\begin{align*} J_n(D) &  =  d^{n}|\alpha_{n,\widehat{x}}|^2\int_D \psi_{n}F_{\widehat{x},n}^{ZZ}\  T\wedge \frac{i}{2}(dZ_{\widehat{x}}\wedge d\overline{Z}_{\widehat{x}}) \\
                &  \ \ \ + 2d^{n}\int_D \psi_{n} \, \mathrm{Re} \left[ \left((q F_{\widehat{x},n}^{ZZ} \alpha_{n,\widehat{x}}\overline{\gamma}_{n,\widehat{x}}\overline{W}_{\widehat{x}}^{q-1}+ F_{\widehat{x},n}^{ZW} \alpha_{n,\widehat{x}}\overline{\beta}_{n,\widehat{x}} \right) \  T\wedge \frac{i}{2}(dZ_{\widehat{x}}\wedge d\overline{W}_{\widehat{x}}) \right] \\
                & \ \ \ \ \ + d^{n}\int_D \psi_{n} \, q^2 F_{\widehat{x},n}^{ZZ} |\gamma_{n,\widehat{x}}|^2|W_{\widehat{x}}|^{2(q-1)} \, T\wedge \frac{i}{2}(dW_{\widehat{x}}\wedge d\overline{W}_{\widehat{x}})  \\
                & \ \ \ \ \ \ \ + d^{n} \int_D 2 \mathrm{Re} \left( q F_{\widehat{x},n}^{ZW} \gamma_{n,\widehat{x}}\overline{\beta}_{n,\widehat{x}}W_{\widehat{x}}^{q-1} \right)  T\wedge \frac{i}{2}(dW_{\widehat{x}}\wedge d\overline{W}_{\widehat{x}}) .
\end{align*}
Let us verify that the first item of Proposition \ref{lemma:decompositionofmuforthesliceconjecture} remains valid in the resonant case. Let $\widehat{x}$ be a $\widehat \mu$-generic element and let $(n_k)_k$ be a sub-sequence such that $\widehat{x}_{-n_k}\in \mathcal{A}_{\tau}$. Assume  $$ \psi_n = \psi\circ f^{-n_k}_{\widehat{x}}\leq b,\ \mu-a.e.\ on\ B $$
and let us verify that there exists $C'(\widehat x)$ such that  for every borel set $U\subset B$:
$$ J_{n_k}(U) \leq C'(\widehat x) e^{-n_k (\lambda_1 + \lambda_2 - \LLog \, d)} e^{2 n_k \varepsilon} . $$
We use the formula on $J_n(D)$ given above and Cauchy-Schwarz inequality. We also use $|F_{\widehat{x},n}^{AB}|\leq \tau$ given by Equation (\ref{eq:equationSIGMA}), item $(i)$ of Lemma \ref{lemma:auxiliarylemma} and $|W_{\widehat{x}}| \leq  \rho_{\varepsilon} := \rho_{\varepsilon}(\hat{x})$ on $B$ provided by Equation \eqref{eq: IB}. If $\sigma_{\widehat{x}} := T\wedge dd^c|W_{\widehat{x}}|^2$, $\lambda_{\widehat{x}} := T\wedge dd^c|Z_{\widehat{x}}|^2$, we indeed obtain
 \begin{align*}
    J_{n_k}(U)  & \leq  \ d^{n_k}|\alpha_{n_k,\widehat{x}}|^2 \,b \tau  \lambda_{\widehat{x}} (U)\\
    & \ \ \ +  \, 2d^{n_k}|\alpha_{n_k,\widehat{x}}| \, b \tau \left( q |\gamma_{n_k,\widehat{x}}| \rho_{\varepsilon}^{q-1} + |\beta_{n_k,\widehat{x}}| \right)  (\lambda_{\widehat{x}} + \sigma_{\widehat{x}})(U)\\
   &  \ \ \ \ \ +  \, qd^{n_k}|\gamma_{n_k,\widehat{x}}| \, b \tau \left( q |\gamma_{n_k,\widehat{x}}| \rho_{\varepsilon}^{2(q-1)} + 2|\beta_{n_k,\widehat{x}}| \rho_{\varepsilon}^{q-1}\right)  \sigma_{\widehat{x}} (U) .
\end{align*}
The existence of $C'(\widehat x)$ then comes from the upper bounds  $|\alpha_{n, \hat{x}}| \leq  e^{- n (\lambda_1 - \varepsilon) }$, $|\beta_{n, \hat{x}}| \leq  e^{-n(\lambda_2 - \varepsilon)}$, $|\gamma_{n, \hat{x}}| \leq M_{\varepsilon}(\widehat{x})  e^{- n (\lambda_1 - \varepsilon)}$ and from $\lambda_1 > \lambda_2$.

\section{Suspensions of one-dimensional Latt\`es maps}\label{sec:relevedeLattes}

Let $\theta=[P:Q]$ be  Latt\`es map on $\mathbb{P}^1$ of degree $d\geq2$ and let $f :=[P:Q:t^d]$ be its suspension on $\mathbb{P}^2$. 
These  maps were studied by Berteloot-Loeb \cite{BL}. We identify the first affine chart with $\mathbb C^2$. 
 Let 
$$G_\theta :=\lim_{n\to+\infty}\frac{1}{d^n} \mathrm{Log}\, || (P^n , Q^n) ||  : \mathbb{C}^2 \to \mathbb R \cup \{ - \infty\} $$
be the Green function of the polynomial mapping $(P,Q)$ and let $G := \max\{G_{\theta},0\}$. The attracting basin $A$ of $(0,0)$ is bounded and equal to $\{ G_\theta < 0 \}$. The Green current of $f$ satisfies $T = dd^c G$ on $\mathbb C^2$, the support of $\mu$ coincides with the boundary $\partial A$ of $A$.

By {\cite[Proposition 3.1]{BL}} (see also \cite{dup03} in higher dimensions), for every $p$ outside a finite number of circles drawn on $\partial A$, there exists a biholomorphism $\mathfrak{p}:(\mathbb{D}^2,0)\longrightarrow (\mathbb{D}^2,p)$ such that $$G_0 (z,w): =G_{\theta}\circ \mathfrak{p}(z,w)=\mathrm{Re}(z) + |w|^2 \textrm{ on } (\mathbb{D}^2,0) .$$
  We denote $$T_0 :=\mathfrak{p}^* T  =dd^c(\max\{G_0,0\}) \textrm { and } \mu_0:=\mathfrak{p}^* \mu  . $$
The closed positive current $T_0$ has the following matrix representation
\begin{equation*}
    T_0=
    \begin{bmatrix}
        T_{11} & 0 \\
        0 & T_{22}
    \end{bmatrix} ,
\end{equation*}
where $T_{11} :=T_0\wedge dd^c|w|^2$ and $T_{22}:=T\wedge dd^c|z|^2$. We describe below $T_0$ and $\mu_0$ on $(\mathbb{D}^2,0)$. We assume to simplify that the germ $(\mathbb{D}^2,0)$ is the open set $D := ]-1,1[^2 \times \mathbb D$. Let 
$$\Omega := \{\mathrm{Re}(z)+|w|^2>0\}  \cap D \textrm { and } M_0:=\{\operatorname{Re}(z)+|w|^2=0\}  \cap D. $$
Observe that $M_0$ is the image (intersected with $D$) of the euclidian $3$-sphere of  $\cmplex^2$ by the classical Cayley transformation, see \cite[Chapter 2.3]{rudin}. Let us parametrize $M_0$ by 
\begin{equation}\label{eq:parametrizationPhi}
    \Phi: ]-1,0]\times ]-1,1[\times ]0,2\pi[ \to M_0 \ \ , \ \ (u,v,\theta) \mapsto (u+iv,\sqrt{-u}e^{i\theta})
\end{equation}
 and let  $\mathrm{Leb}_{M_0}:=\Phi_*\mathrm{Leb}_{\reels^3}$. Let also $\omega_0$ be the standard hermitian $(1,1)-$form of $\cmplex^2$.

\begin{prop}\label{prop:featureofTnearregularpoints} With the preceding notations, 
$$T_0=\begin{bmatrix}
\frac{1}{8}\mathrm{Leb}_{M_0} & 0 \\
0 & \mathrm{Leb}_{\Omega} + \frac{|w|^2}{2} \mathrm{Leb}_{M_0}
\end{bmatrix},$$
and 
$$\mu_0=\frac{1}{8}\mathrm{Leb}_{M_0}=T_0\wedge dd^c\left|w\right|^2 . $$
In particular, $$\mu_0=\psi_0 \left(T_0\wedge\omega_0\right) , \textrm { where } \psi_0(z,w):=\frac{1}{1+4 |w|^2}\displaystyle{1\!\!1}_{M_0} . $$ 
\end{prop}

The remainder of this Section is devoted to the proof of Proposition {\ref{prop:featureofTnearregularpoints}}. \\

\textit{Computation of $T_{11}$.--} Let $\varphi$ be a test function with compact support in $D$. We denote $\varphi_{z\overline{z}} := \frac{\partial}{\partial z}\frac{\partial}{\partial \overline{z}}\varphi$, similarly for $\varphi_{w\overline{w}}$. By definition we have:
\begin{align*}
    \langle T_{11}, \varphi\rangle & = \int_{\Omega} G_0(z,w)\varphi_{z\overline{z}}(z,w)\ \mathrm{d}\mathrm{Leb}(z,w).
\end{align*}
Let $\Delta_w$ denote the subset $]-1,1[ \times ]-1,1[ \times\{w\}$ intersected with $\Omega$, so that
\begin{equation}\label{eq:crochetdeT_11}
    \langle T_{11}, \varphi\rangle  = \int_{\mathbb{D}} \left [ \int_{\Delta_w} G_0(z,w)\varphi_{z\overline{z}}(z,w)\ \mathrm{d}\mathrm{Leb}(z) \right] \mathrm{d}\mathrm{Leb}(w) =: \int_{\mathbb{D}} I_w(\varphi)\ \mathrm{d}\mathrm{Leb}(w).
\end{equation}
We write $z=u+iv$, hence $\varphi_{z\overline{z}} = \frac{1}{4} \Delta_{u,v} \varphi$. Observe also that $G_0$ only depends on $(u,w)$. The number $4 I_w(\varphi)$ is then equal to
$$ \int_{-1}^{1} \left[ \int_{-|w|^2}^{1}G_0(u,w)\frac{\partial^2 \varphi}{\partial u^2}(u+iv,w)\ \mathrm{d}u \right] \mathrm{d}v + \int_{-|w|^2}^{1}G_0(u,w) \left[ \int_{-1}^{1}\frac{\partial^2 \varphi}{\partial v^2}(u+iv,w)\ \mathrm{d}v \right] \mathrm{d}u. $$
The second term vanishes since $\mathrm{Supp}(\varphi(\cdot,w))\subset ]-1,1[ \times ]-1,1[$. Integrating by parts and using $\partial_u{G_0}=1$, $\partial_v G_0=0$ for the first  one, we get
\begin{align*}
   4 I_w(\varphi) =  -  \int_{-1}^{1} \left[ \int_{-|w|^2}^{1} \frac{\partial \varphi}{\partial u}(u+iv,w)\ \mathrm{d}u\right] \mathrm{d}v\ .
\end{align*}
Using again $\mathrm{Supp}(\varphi(\cdot,w))\subset ]-1,1[ \times ]-1,1[$, we get
\begin{align*}
    I_w(\varphi) = \frac{1}{4} \int_{-1}^{1}\varphi(-|w|^2+iv,w)\ \mathrm{d}v.
\end{align*}
We therefore obtain by  (\ref{eq:crochetdeT_11}): 
\begin{align*}
    \langle T_{11}, \varphi\rangle & = \frac{1}{4} \int_{\mathbb{D}} \int_{-1}^{1}\varphi(-|w|^2+iv,w)\ \mathrm{d}v\ \mathrm{d}\mathrm{Leb}(w) 
 \end{align*}
Observe that $(w,v)\mapsto (-|w|^2+iv,w)$ is another parametrization of $M_0$ different from $\Phi$ defined in (\ref{eq:parametrizationPhi}). Let us write this integral with polar coordinates:
\begin{align*}
    \langle T_{11}, \varphi\rangle & = \frac{1}{8}\int_{-1}^{1}\int_{0}^{2\pi}\int_{0}^{1}\varphi(-\rho^2+iv,\rho e^{i\theta})\ 2\rho\mathrm{
    d}\rho\ \mathrm{d}\theta\ \mathrm{d}v.
\end{align*}
Applying the change of variable $u=-\rho^2$ one gets back to the parametrization $\Phi$:
$$    \langle T_{11}, \varphi\rangle  = \frac{1}{8}\int_{]-1,0]\times ]-1,1[\times ]0,2\pi[}\varphi(u+iv,\sqrt{-u}e^{i\theta})\ \mathrm{d}\mathrm{Leb}(u,v,\theta) = \frac{1}{8}\int_{M_0} \varphi\ \mathrm{d}\mathrm{Leb}_{M_0} . $$

\textit{Computation of $T_{22}$.--} Let $\Delta_z$ be the vertical disc $\{z\}\times\mathbb{D}$ intersected with $\Omega$. One has:
$$   \langle T_{22}, \varphi\rangle = \int_{]-1,1[ \times ]-1,1[} \left[ \int_{\Delta_z }G_0(z,w)\varphi_{w\overline{w}}(z,w)  \mathrm{d}\mathrm{Leb}(w) \right] \mathrm{d}\mathrm{Leb}(z) . $$
Writing  $z=u+iv$, we get
$$  \langle T_{22}, \varphi\rangle = \int_{-1}^{1} \left[ \int_{-1}^{0}J_{u,v}(\varphi)\ \mathrm{d}u \right]  \mathrm{d}v + \int_{-1}^{1} \left[ \int_{0}^{1}J_{u,v}(\varphi)\ \mathrm{d}u \right] \mathrm{d}v, $$
where 
$$J_{u,v}(\varphi):=\int_{\mathbb{D}\backslash\mathbb{D}(\sqrt{-u})}G_0(u+iv,w)\varphi_{w\overline{w}}(u+iv,w)\ \mathrm{d}\mathrm{Leb}(w) ,$$
with the convention $\mathbb{D}(\sqrt{-u}) = \emptyset$ if $0 \leq u \leq 1$. 

\begin{lemme}\label{coupe} $ $
\begin{enumerate}
\item $J_{u,v}(\varphi) = \int_{\mathbb{D}}\varphi(u+iv,w)\ \mathrm{d}\mathrm{Leb}(w)$ for $0 \leq u \leq 1$. 
\item $J_{u,v}(\varphi) = \int_{\mathbb{D}\backslash\mathbb{D}(\sqrt{-u})}\varphi(z,w)\ \mathrm{d}\mathrm{Leb}(w) -\frac{u}{2}\int_{0}^{2\pi}\varphi(z,\sqrt{-u}e^{i\theta})\ \mathrm{d}\theta$ for $-1 \leq u \leq 0$.
\end{enumerate}
\end{lemme}

\noindent\textbf{\underline{Proof:}} We write $w=s+it$ and use  $\varphi_{w\overline{w}}=\frac{1}{4}  \Delta_{s,t} \varphi$. Since $\varphi$ has compact support in $D$, one gets by Green-Stokes formula :
$$ 4 J_{u,v}(\varphi):=\int_{\mathbb{D}} G_0(u+iv,w)  \Delta_{s,t} \varphi \ \mathrm{d}\mathrm{Leb}(w) = \int_{\mathbb{D}} \Delta_{s,t} G_0(u+iv,w) \varphi \ \mathrm{d}\mathrm{Leb}(w)$$
The first item follows since $\partial^2_sG_0=\partial_t^2G_0=2$. In the case $-1 \leq u \leq 0$, we have:
\begin{equation}\label{eq: rappel}
    J_{u,v}(\varphi) = \int_{\mathbb{D}\backslash\mathbb{D}(\sqrt{-u})} G_0(u+iv,w)\varphi_{w\overline{w}}(u+iv,w)\ \mathrm{d}\mathrm{Leb}(w).
\end{equation}
Let us denote $z=u+iv$ and $w=s+it\in \mathbb{D}\backslash \mathbb{D}(\sqrt{-u})$. Using $\mathrm{Supp}(\varphi(z,\cdot))\subset\mathbb{D}$ and $G_0(z,\cdot)|_{\partial\mathbb{D}(\sqrt{-u})}=0$, Green-Stokes formula gives :
\begin{align*}
    & \int_{\mathbb{D}\backslash\mathbb{D}(\sqrt{-u})}  G_0(z,s+it)\Delta_{s,t}\varphi(z,s+it) - \Delta_{s,t}G_0(z,s+it) \varphi(z, s+it) \ \mathrm{d}\mathrm{Leb}(s,t) \\
                                & = \int_{\partial\mathbb{D}(\sqrt{-u})} \varphi(z,s+it)\left(\nabla_{(s,t)}G_0(z,s+it)\right)\cdot (s,t)\ \mathrm{d}\sigma(s,t) ,
\end{align*}
where $\sigma$ is the Lebesgue measure on the circle $\partial\mathbb{D}(\sqrt{-u})$. Since $\Delta_{s,t}G_0(z,s+it)=4$ and $\left(\nabla_{(s,t)}G_0(z,s+it)\right)\cdot (s,t)=2|w|^2=-2u$, we infer 
\begin{align*}
    4 \int_{\mathbb{D}\backslash\mathbb{D}(\sqrt{-u})} \varphi(z,w) -G_0(z,w)\varphi_{w\overline{w}}(z,w) \ \mathrm{d}\mathrm{Leb}(w) = \int_{\partial\mathbb{D}(\sqrt{-u})}2u\varphi(z,w)\ \mathrm{d}\sigma(w) .
\end{align*}
The second item follows by combining this formula with \eqref{eq: rappel}. \qed \\

Let us now return to the computation of $\langle T_{22}, \varphi\rangle$. Let $\Omega^+:=\Omega\cap \{u>0\}$ and  $\Omega^{-}:=\Omega\cap\{u<0\}$. According to Lemma \ref{coupe}, we get 
$$  \int_{-1}^{1}\int_{0}^{1}J_{u,v}(\varphi)\ \mathrm{d}u\ \mathrm{d}v = \int_{\Omega^+}\varphi(z,w)\ \mathrm{d}\mathrm{Leb}(z,w) , $$
and that $ \int_{-1}^{1}\int_{-1}^{0}J_{u,v}(\varphi)\ \mathrm{d}u\ \mathrm{d}v $ is equal to 
$$ \int_{\Omega^{-}}\varphi(z,w)\ \mathrm{d}\mathrm{Leb}(z,w)\\
                 + \int_{-1}^{1} \left[ \int_{-1}^{0} \left[ \int_{0}^{2\pi}-\frac{u}{2}\varphi(u+iv,\sqrt{-u}e^{-i\theta})\ \mathrm{d}\theta \right] \mathrm{d}u \right] \mathrm{d}v . $$
If $h(z,w):=\frac{|w|^2}{2}$ the second integral is equal to:
$$\int_{]-1,0[\times]-1,1[\times ]0,2\pi[} (h \varphi) \circ \Phi(u,v,\theta) \ \mathrm{d}(u,v,\theta) = \int_{M_0}h\varphi\ \mathrm{d}\mathrm{Leb}_{M_0} . $$
We deduce as desired
\begin{align*}
    \langle T_{22}, \varphi\rangle = \int_{-1}^{1}\int_{-1}^{1} J_{u,v}(\varphi)\ \mathrm{d}u\ \mathrm{d}v = \int_{\Omega}\varphi(z,w)\ \mathrm{d}\mathrm{Leb}(z,w) + \int_{M_0}h\varphi\ \mathrm{d}\mathrm{Leb}_{M_0}.
\end{align*}

\textit{Computation of $\mu_0$.--} Writing $G_0^+:=\max\{G_0,0\}$, we have $T_0=dd^c G_0^+$ and $\mu_0=dd^c(G_0^+dd^cG_0^+)$. Hence
\begin{align*}
    \langle \mu_0,\varphi \rangle & = \langle T_0,G_0^+dd^c\varphi \rangle = \langle T_{11}dd^c|z|^2 + T_{22}dd^c|w|^2,G_0^+dd^c\varphi\rangle\\
                                  & = \frac{1}{8}\int_{M_0}G_0\varphi_{w\overline{w}}\ \mathrm{d}\mathrm{Leb}_{M_0} + \int_{\Omega}G_0\varphi_{z\overline{z}}\ \mathrm{d}\mathrm{Leb} +  \int_{M_0}G_0\varphi_{z\overline{z}}h\ \mathrm{d}\mathrm{Leb}_{M_0}.
\end{align*}
Since $G_0=0$ on $M_0$, we get $ \langle \mu_0,\varphi \rangle = \int_{\Omega}G_0\varphi_{z\overline{z}}\ \mathrm{d}\mathrm{Leb}_{\Omega} = \langle T_{11}, \varphi\rangle$ as desired (observe that this is the Lebesgue part on $\Omega$ in $T_{22}$ which provides  $\mu = T_{11}$). The formula $\mu_0=\frac{1}{1+8h}\displaystyle{1\!\!1}_{M_0} T_0\wedge\omega_0$ immediately follows.

\begin{otherlanguage}{english}
\bibliographystyle{abbrv}
\bibliography{biblioPROC} 
\end{otherlanguage}

$ $ \\
\noindent {\footnotesize C. Dupont, V. Tapiero}\\
{\footnotesize Universit\'e de Rennes}\\
{\footnotesize CNRS, IRMAR - UMR 6625}\\
{\footnotesize 35042 Campus de Beaulieu, Rennes, France}\\
{\footnotesize christophe.dupont@univ-rennes.fr}\\
{\footnotesize virgile.tapiero@univ-rennes.fr}

\end{document}